\documentclass[10pt]{amsart}
\usepackage{amssymb,amsmath,amscd,amsfonts,amsthm,mathrsfs}
\usepackage{verbatim} 
\usepackage{a4wide}
\usepackage{color}

\usepackage[latin1]{inputenc}

\newcommand{\N}{\mathbb{N}}

\newcommand{\R}{\mathbb{R}}
\newcommand{\C}{\mathbb{C}}
\newcommand{\Bc}{\mathcal{B}}
\newcommand{\Cc}{\mathcal{C}}
\newcommand{\Dc}{\mathcal{D}}
\newcommand{\Ec}{\mathcal{E}}
\newcommand{\Hc}{\mathcal{H}}
\newcommand{\Ic}{\mathcal{I}}
\newcommand{\Jc}{\mathcal{J}}
\newcommand{\Nc}{\mathcal{N}}
\newcommand{\Oc}{\mathcal{O}}
\newcommand{\Id}{{\rm{Id}} }

\newcommand{\Sc}{\mathcal{S}}
\def\cchi{\raisebox{.45 ex}{$\chi$}}
\newcommand{\ad}{\mathrm{ad}}  
\newcommand{\supp}{\mathrm{supp}}

\newcommand{\Dr}{\mathscr{D}}

\newcommand{\Gr}{\mathscr{G}}
\newcommand{\Hr}{\mathscr{H}}
\newcommand{\Kr}{\mathscr{K}}
\newcommand{\Sr}{\mathscr{S}}

\newtheorem{theorem}{Theorem}[section]
\newtheorem{proposition}{Proposition}[section]
\newtheorem{lemma}{Lemma}[section]
\newtheorem{remark}{Remark}[section]
\newtheorem{example}{Example}[section]

\numberwithin{equation}{section}

\def \e{\varepsilon}

\def \bone{\mathbf{1}}

\def \rmi{{\rm i}}

\title[LAP and propagation of long range Dirac systems]{Limiting 
absorption principle for some long range 
perturbations of Dirac systems at threshold energies}

\begin{document}
\author{Nabile Boussaid}
\address{D\'epartement de Math\'ematiques,
Universit\'e de Franche-Comt\'e,
16 route de Gray,
25030 Besan\c{c}on Cedex,
France}
\email{nabile.bousssaid@univ-fcomte.fr}
\author{Sylvain Gol\'enia}
\address{Mathematisches Institut der Universit\"at Erlangen-N\"urnberg,
Bismarckstr.\ 1 1/2 \\
91054 Erlangen, Germany}
\email{golenia@mi.uni-erlangen.de}
\date{Version of \today}
\begin{abstract}
We establish a limiting absorption principle for some long range 
perturbations of the Dirac systems at threshold energies. We cover
multi-center interactions with small coupling constants. 
The analysis is reduced to study a family of non-self-adjoint
operators. The technique is based on a positive commutator theory
for non self-adjoint operators, which we develop in appendix. We also discuss some
applications to the dispersive Helmholzt model in the quantum regime.  
\end{abstract}
\maketitle
\tableofcontents
\section{Introduction}
We study properties of relativistic massive charged particles with
spin-$1/2$ (e.g., electron, positron, (anti-)muon,
(anti-)tauon,$\ldots$). We follow the Dirac formalism, see \cite{Dirac}. Because of the spin, the configuration space of the
particle is vector valued. To simplify, we consider finite
dimensional and trivial fiber. Let $\nu\geq 2$ be an integer. 
The movement of the free particle is given by the Dirac equation, 
\begin{equation*}
i \hbar \frac{\partial \varphi}{\partial t }= D_m \varphi, \mbox{ in } L^2(\R^3; \C^{2\nu}),
\end{equation*}
where $m>0$ is the mass, $c$ the speed of light, $\hbar$ the
reduced Planck constant, and 
\begin{equation}\label{e:op}
D_m:= c \hbar\, \alpha\cdot P+mc^2\beta=
-\rmi c\hbar\sum_{k=1}^3\alpha_k\partial_k + mc^2\beta. 
\end{equation}
Here we set ${\alpha}:=\left(\alpha_1, \alpha_2,  \alpha_3\right)$
and $\beta:=\alpha_4$.  The  $\alpha_i$, for $i\in\{1,2,3, 4\}$,
are linearly independent self-adjoint linear applications, acting in $\C^{2\nu}$, satisfying the anti-commutation relations: 
\begin{equation}\label{e:Diracrep0}
\alpha_i\alpha_j+\alpha_j\alpha_i=2\delta_{ij}\bone_{\C^{2\nu}}, \mbox{
  where } i,j\in\{1,2,3,4\}.
\end{equation}
For instance, when $\nu=2$, one may choose the Pauli-Dirac
representation: 
\begin{gather}\label{e:Diracrep}
\alpha_i=\left(
\begin{array}{cc}
0&\sigma_i\\
\sigma_i&0
\end{array} \right)\quad\mbox{ and }\quad
\beta=\left(\begin{array}{cc}
\Id_{\C^\nu}&0\\
0&-\Id_{\C^\nu}
\end{array} \right)
\\ \nonumber
\mbox{ where  } \sigma_1=\left(
\begin{array}{cc}
0&\; 1\\
1&\;0
\end{array} \right),\quad \quad\sigma_2=\left(
\begin{array}{cc}
0&-\rmi\\
\rmi&0
\end{array} \right)\quad\mbox{ and } \quad\sigma_3=\left(
\begin{array}{cc}
1&0\\
0&-1
\end{array} \right),
\end{gather}
for $i=1,2,3$. We refer to \cite{Thaller}[Appendix 1.A] for
various equivalent representations. In this paper we do not choose
any specific basis and work intrinsically with
\eqref{e:Diracrep0}. We refer to \cite{LawsonMichelsohn} for a
discussion of the representations of the Clifford algebra generated by
\eqref{e:Diracrep0}. We also renormalize and consider $\hbar=c=1$. 
The operator $D_m$ is essentially self-adjoint on $\Cc^\infty_c(\R^3;
\C^{2\nu})$ and the domain of its closure is $\Hr^1(\R^3; \C^{2\nu})$,
the Sobolev space of order $1$ with values in $\C^{2\nu}$. We denote
the closure with the same symbol. Easily, using Fourier
transformation and some symmetries, one deduces the spectrum of $D_m$ 
is purely absolutely continuous and given by $(-\infty, -m]\cup[m,
    \infty)$.  

In this introduction, we focus on the dynamical and spectral properties of
the Hamiltonian describing the movement of the particle interacting
with $n$ fixed, charged particles. We model them by fixed points
$\{a_i\}_{i=1,\ldots, n}\in \R^{3n}$ with respective charges 
$\{z_i\}_{i=1,\ldots, n}\in \R^{n}$. Doing so, we tacitly suppose that
the particles $\{a_i\}$ are far enough from one another, so as to
neglect their interaction.  Note we make no hypothesis on the sign of the
charges. The new Hamiltonian is given by   
\begin{eqnarray}\label{e:H}
H_\gamma:=D_m + \gamma V_c(Q), \mbox{ where } V_c:=v_c
\otimes\Id_{\C^{2\nu}} \mbox{ and } v_c(x):=\sum_{k=1,\ldots ,n}
\frac{z_i}{|x-a_i|},
\end{eqnarray} 
acting on $\mathcal{C}_c^\infty(\R^3\setminus\{a_i\}_{i=1, \ldots, n};
\C^{2\nu})$, with $a_i\neq a_j$ for $i\neq j$. The $\gamma\in \R$ is
the coupling constant. The index $c$  stands for \emph{coulombic
 multi-center}. The notation $V(Q)$ indicates 
the operator of multiplication by $V$. Here, we identify
$L^2(\R^3;\C^{2\nu})\simeq L^2(\R^3)\otimes \C^{2\nu}$,
canonically. Remark the perturbation $V_c$ is not relatively compact
with respect to $D_m$, then one needs to be careful to define a
self-adjoint extension for $D_m$. Assuming
\begin{eqnarray}\label{e:atomic} 
Z:=|\gamma| \max_{i=1, \ldots, n}(|z_i|)< \sqrt{3}/2,
\end{eqnarray} 
the theorem of Levitan-Otelbaev ensures that 
$H_\gamma$ is essentially self-adjoint and its domain is
the Sobolev space $\Hr^1(\R^3; \C^{2\nu})$, see \cite{AraiYamada,
  Kalf,  Klaus, LandgrenRejto, LandgrenRejtoKlaus, LevitanOtelbaev}
for various generalizations. This condition corresponds to the 
 nuclear charge $\alpha_{\rm at}^{-1} Z\leq 118 $,
where $\alpha_{\rm at}^{-1}= 137.035999710(96)$.  Note that using the
Hardy-inequality, the Kato-Rellich theorem will apply till $Z<1/2$
and is optimal in the matrix-valued case, see \cite{Thaller}[Section
  4.3] for instance. For $Z<1$, one 
shows there exists only one self-adjoint 
extension so that its domain is included in $\Hr^{1/2}(\R^3;
\C^{2\nu})$, see \cite{Nenciu}. This covers the nuclear charges up to
$Z=137$. When $n=1$ and $Z=1$, this property still
holds true, see \cite{EstebanLoss}. Surprisingly enough, when $n=1$
and $Z>1$, there is no self-adjoint extension with domain included
in $\Hr^{1/2}(\R^3; \C^{2\nu})$, see \cite{Xia}[Theorem 6.3]. We mention also the
work of \cite{VoronovGitmanTyutin} for $Z>1$. 

In \cite{Nenciu} and $Z<1$, one shows the essential
spectrum is given by  $(-\infty, -m]\cap [m, \infty)$ for all
self-adjoint extension. For all $Z$, one refers to
\cite{GeorgescuMantoiu}[Proposition 4.8.], which relies on
\cite{Xia}. In \cite{GoleniaGeorgescu}  one gives some criteria of
stability of the essential spectrum for 
some very singular cases. In \cite{BerthierGeorgescu}, one
 proves there is no embedded 
eigenvalues for a more general model and till the coupling constant
$Z<1$. For all energies being in a
compact set included in $(-\infty, -m)\cap (m, \infty)$, 
\cite{GeorgescuMantoiu} obtains some estimates of the resolvent. 
This implies some propagation estimates and that the
spectrum of $H_\gamma$ is purely absolutely continuous. Similar
results have been obtained  for magnetic potential of constant
direction, see \cite{Yokoyama}
and more recently \cite{RichardTiedra}.

In this paper we are interested in  uniform estimates of the resolvent
at threshold energies. The energy $m$ is called the \emph{electronic threshold}
and $-m$ the \emph{positronic threshold}. 
In Theorem \eqref{t:main}, we obtain a uniform estimation of the 
resolvent over $[-m-\delta, -m]\cup [m, m+\delta]$, see \eqref{e:main}
and deduce some propagation properties, see \eqref{e:Kato}. One
difficulty is that in the case $n=1$ and $z_i<0$, 
it is well known there are infinitely many eigenvalues  in the gap
$(-m, m)$ converging to the $m$ as soon as $\gamma \neq 0$ (see 
for instance \cite{Thaller}[Section 7.4] and references therein).  
This is a difficult problem and, to our knowledge, this result is new
for the multi-center case.  There is a larger literature
for non-relativist models, e.g., $-\Delta +V$ in $L^2(\R^n;\C)$. The
question is intimately linked with the presence of resonances at
threshold energy, \cite{JensenNenciu, FournaisSkibsted, Nakamura,
  Richard, Yafaev82}. We mention also
\cite{BurqPlanchonStalkerTahvildarZadeh} for applications to
Strichartz estimates and \cite{DerezinskiSkibsted} for
  applications to scattering theory. 

Before giving the main result, we shall discuss some commutator methods. The 
first stone was set C.R.\ Putnam a self-adjoint operator $H$
acting in a Hilbert space $\Hr$, see \cite{Putnam} and for instance
\cite{ReedSimon}[Theorem XIII.28]. One supposes there is a \emph{bounded}
operator $A$ so
that  \begin{equation}\label{e:KatoPutnam} 
C:=[H,iA]_\circ>0,
\end{equation}
where $>$ means non-negative and injective. The commutator has to be
understood in the form sense. When it extends into a bounded operator
between some spaces, we denote this extension with the symbol $\circ$
in subscript, see Appendix \ref{s:dev-commut}. The operator 
$A$ is said to be \emph{conjugate} to $H$. One deduces some estimation
on the imaginary part of the resolvent, i.e., one finds some
\emph{weight} $B$, a closed injective operator with dense domain, so that 
\begin{eqnarray*}
\sup_{\Re(z)\in \R, \Im(z)>0} \Im \langle f, (H- z )^{-1} f \rangle \leq \|B f\|^2.
\end{eqnarray*}
This estimation is equivalent to the global propagation estimate,
c.f.\ \cite{Kato} and \cite{ReedSimon}[Theorem XIII.25]: 
\begin{eqnarray*}
\int_\R \|B^{-1} e^{itH} f\|^2 dt \leq 2 \|f\|^2
\end{eqnarray*}
One infers that the spectrum of $H$ is purely absolutely
continuous with respect to the Lebesgue measure. In particular, $H$
has no eigenvalue.  To deal with the presence of eigenvalues, the fact
that $A$ is bounded and with the  $3$-body-problem, Eric Mourre has
the idea to localized in energy the estimates and to allow a compact
perturbation, see \cite{Mourre}.  With further hypothesis, one shows
an estimate of the resolvent (and not only on the imaginary
part). The applications of this theory are numerous. The theory was
immediately adapted to treat the $N$-body problem, see
\cite{PerrySigalSimon}. The theory was finally improved in many 
directions and optimized in many ways, see
\cite{AmreinBoutetdeMonvelGeorgescu} for a more thorough discussion of
these matters. We mention also \cite{GeorgescuGerardMoller,
  GoleniaJecko, Gerard} for recent 
developments. As we are concern about thresholds, Mourre's method
does not seem enough, as the estimate of the resolvent is given on an
interval which is strictly smaller than the one used in the commutator
estimate. In \cite{BoutetMantoiu} one generalizes the result of
Kato-Putnam approach. Under some conditions, one allows $A$ to be
unbounded. They obtain a global estimate of the resolvent. Note this 
implies the absence of eigenvalue. In \cite{FournaisSkibsted}, in the
non-relativistic context, by asking some positivity on the Virial of
the potential, see below, one is able to conciliate the estimation of
the resolvent above the threshold energy and the accumulation of eigenvalue
under it. In \cite{Richard}, one presents an abstract version of the method of
\cite{FournaisSkibsted}.  To give an idea, we shall
compare the theories on a non-optimal example. Take
$H:=-\Delta+V$ in $L^2(\R^3)$, with $V$ being in the Schwartz
space. Consider the generator of dilation $A:= (P\cdot Q + Q\cdot
P)/2$, one looks at the quantity 
\begin{eqnarray*}
[H, iA]_\circ - cH = -(2-c) \Delta - W_V(Q), \mbox{ where } W_V(Q):= Q \cdot
\nabla V(Q) + c V(Q),  
\end{eqnarray*} 
with $c\in (0,2)$ and seeks some positivity. The expression $W_V$ is
called the \emph{Virial} of $V$. In \cite{FournaisSkibsted}, one uses
extensively that $W_V(x)\leq  -c\langle x\rangle^{-\alpha}$ for some
$\alpha,c>0$ and $|x|$ big enough. In \cite{Richard}, one notices that
it suffices to suppose that $W_V(x)\leq 0$ and to take advantage
of the positivity of the Laplacian. We take the opportunity to mention
that it is enough to suppose that $W_V(x)\leq c'|x|^{-2}$, for some
small positive constant $c'$, see Theorem
\ref{t:nonrelamain}. Observe also that these methods give different 
weights. For instance, \cite{FournaisSkibsted} obtains better weights
in the scale of $\langle Q\rangle^\alpha$ and \cite{Richard} can
obtain singular weights like $|Q|$, see Appendix
\ref{s:WeakMourreTheory}. Finally, \cite{FournaisSkibsted} 
deals only with low energy estimates and \cite{Richard} works globally
on $[0, \infty)$. We also point  out \cite{Herbst} which relies on
commutator techniques and deals with smooth homogeneous potentials. 

In this article, we revisit the approach of \cite{Richard} and make
several improvements, see Appendix \ref{s:WeakMourreTheory}. 
Our aim is twofold: to treat dispersive non self-adjoint operator and
to obtain estimates of the resolvent uniformly in a parameter.  
At first sight, these improvements are pointless from the standpoint
of the Coulomb-Dirac problem we treat. In reality, they are the
key-stone of our approach. 

As a direct by-product of the method, we obtain some new results for
dispersive Schr\"odinger operators. The following $V_2$ term corresponds to the
absorption coefficient of the laser energy by material medium
absorption term in the Helmholtz model, see \cite{Jackson} for  
instance. 

\begin{theorem}\label{t:nonrela}
Let $n\geq 3$. Suppose that $V_1, V_2\in L^\infty(\R^n;\R)$ satisfy:
\begin{enumerate}
\item[(H1)] $\nabla V_i$,  $ Q\cdot \nabla V_i(Q)$,
 $\langle Q\rangle (Q\cdot \nabla V_i)^2(Q)$ are bounded, for $i\in\{1,2\}$. 
\item[(H2)] There are $c_1\in [0,2)$ and $\displaystyle c_1'\in
  \big[0, 4(2-c_1)/(n-2)^2\big)$ such that 
\begin{eqnarray*}
W_{1}(x):= x\cdot
  (\nabla V_1)(x) + c_{1} V_1(x) \leq \frac{c_{1}'}{|x|^2}, \mbox{ for
  all } x\in \R^n. 
\end{eqnarray*} 
and 
\begin{eqnarray*}
V_2(x)\geq 0 \mbox{ and } -x\cdot
  (\nabla V_2)(x)\geq  0, \mbox{ for all } x\in \R^n.
\end{eqnarray*} 
\end{enumerate}
On $\Cc^\infty_c(\R^n)$, we define $H:= -\Delta + V(Q), \mbox{ where }
V:=V_1 + i V_2$. The closure of $H$ defines a dispersive closed operator
with domain $\Hr^2(\R^n)$. We keep denoting it with $H$. Its spectrum
included in the upper half-plane. The operator $H$ has no eigenvalue
in $[0,\infty)$. Moreover,  
\begin{eqnarray}\label{e:nonrela}
\sup_{\lambda\in [0,\infty), \, \mu>0} \big\|\, |Q|^{-1}(H -\lambda +
  \rmi \mu )^{-1} |Q|^{-1}\big\| <\infty.	  
\end{eqnarray}	
\end{theorem}
Note we require nor smoothness on the potentials neither that they are
relatively compact with respect to the Laplacian. We refer to  Appendix
\ref{s:nonrela} for further comments, the case $c_1=0$ and a stronger result. 

We come back to the main application, namely the operator $H_\gamma$
defined by \eqref{e:H}. As the Dirac operator is vector-valued,
coulombic  interaction are singular and as we are interested in both
thresholds, we were not able to use directly the ideas of
\cite{FournaisSkibsted, 
Richard}. Indeed, it is unclear for us if one can actually 
deal with thresholds energy  and keep the ``positivity'' of the something close to the quantity $[H_\gamma, iA]- c H_\gamma$, for some self-adjoint operator $A$. We hedge this fundamental problem. 
First of all we cut-off the singularities of the potential $V_c$ and
consider the operator $H^{\rm bd}_\gamma = D_m+\gamma V$ in Section \ref{s:red}.
We recover the singularities of the operator by perturbation in
Proposition \ref{p:collage}. In Section \eqref{s:reduction}, we
explicit the resolvent of  $H^{\rm bd}_\gamma-z$ relatively to a spin-down/up
decomposition.  This transfers the analysis to the one of
an elliptic operator of second order, $\Delta_{m,v,z}$, see Section
\ref{s:reduction}. The drawback is  that this operator is dispersive
and also depends on the spectral parameter $z$. We elude the latter by studying the family $\{\Delta_{m,v,\xi}\}_{\xi\in\Ec}$ uniformly
in $\Ec$.  In Section \ref{s:LAPano}, we explain how
to deduce the estimation of the resolvent of $H^{\rm bd}_\gamma$
having the one of $\Delta_{m,v,z}$. In the Section \ref{s:positive},
we establish some positive commutator estimates for $\Delta_{m,v,z}$
and derive the sought estimates of the resolvent, see Theorem
\eqref{t:mourreD}. For the last step, we rely on the theory developed
in Appendix \ref{s:WeakMourreTheory}. The main result of this
introduction is the following one.   

\begin{theorem}\label{t:main} 
There are $\kappa, \delta, C>0$ such that
\begin{eqnarray}\label{e:main}
\sup_{|\lambda|\in [m, m+\delta], \, \varepsilon>0, |\gamma|\leq \kappa }\|\langle
Q\rangle^{-1}(H_\gamma -\lambda - \rmi \varepsilon )^{-1} \langle
Q\rangle^{-1}\|\leq C.	 
\end{eqnarray}	
In particular, $H_\gamma$ has no eigenvalue in $\pm m$. Moreover,
there is $C'$ so that 
\begin{eqnarray}\label{e:Kato}
\sup_{|\gamma|\leq \kappa}\int_\R \| \langle Q \rangle^{-1}
e^{-it H_\gamma } E_\Ic(H_\gamma) f \|^2 dt  \leq C'\|f\|^2,
\end{eqnarray} 
where $\Ic=[-m-\delta, -m]\cup [m, m+\delta]$ and where
$E_\Ic(H_\gamma)$ denotes the spectral measure of $H_\gamma$. 
\end{theorem}

A more general result is given in Theorem \ref{t:vraimain}. 
In Theorem \ref{t:vraimain2} we discuss the weights $\langle P\rangle^{1/2}
|Q|$ and in Remark \ref{r:noLAP} the weights $|Q|$. If one
is not interested in the uniformity in the coupling constant, 
using \cite{GeorgescuMantoiu}, one can consider all $\delta>0$ and
deduce \eqref{e:main}. The propagation estimate \eqref{e:Kato} 
refers as Kato smoothness and it is a well-known consequence of
\eqref{e:main}, see \cite{Kato}. 
Using some kernel estimates, one can obtain \eqref{e:main} directly
for the free Dirac operator, i.e., $\gamma=0$, see for instance
\cite{Thaller}[Section 1.E] and  \cite{KatoYajima}. One may find an 
alternative proof of this fact in \cite{IftimoviciMantoiu} which
relies on some positive commutator techniques. 

In this study, we are mainly interested by long range perturbations of
Dirac operators. Concerning limiting absorption principle for short
range perturbations of Dirac operators there are some interesting
works such as \cite{DAnconaFanelli} for small perturbations without
discrete spectrum or \cite{Boussaid} for potentials producing discrete
spectrum. These authors were mainly interested by time decay estimates
similar to \eqref{e:Kato}. In the short range case, the limiting
absorption principle is a key ingredient to establish Strichartz
estimates for perturbed Dirac type equations see \cite{Boussaid2,
  DAnconaFanelli2}. For free Dirac equations there are some direct
proofs, see \cite{EscobedoVega, MachiharaNakamuraOzawa,
  MachiharaNakamuraNakanishiOzawa}. Time decay estimates such as
\eqref{e:Kato} or Strichartz are crucial tool to establish well
posedness results \cite{EscobedoVega, MachiharaNakamuraOzawa,
  MachiharaNakamuraNakanishiOzawa} and stability results
\cite{Boussaid, Boussaid2} for nonlinear Dirac equations.

The paper is organized as follows. In the second section we reduced
the analysis of the resolvent of the Dirac operator perturbed with a
bounded potential to the one of family of non self-adjoint
operators. In the third part, we analyze these operators and obtain
some estimates of the resolvent. In the fourth part, we state the main
results of the paper. For the convenience of the reader, we expose
some commutator expansions in the Appendix \ref{s:dev-commut}. In the
Appendix \ref{s:WeakMourreTheory}, we develop the abstract positive
commutator theory. At last in Appendix \ref{s:nonrela}, we give a
direct application to the theory in the context of the Helmholtz
equation.

\noindent
{\bf Notation:} In the following $\Re$ and $\Im$ denote the real and 
imaginary part, respectively. The smooth function with compact support
are denoted by $\Cc^\infty_c$. 
Given a complex-valued function $F$, we denote
by $F(Q)$ the operator of multiplication by $F$. We mention also the
notation $P=-\rmi \nabla$. We use the standard $\langle \cdot
\rangle:= (1+ |\cdot|^2)^{1/2}$.   

\noindent
{\bf Acknowledgments: } We would like to thank Lyonel Boulton,
Bertfried Fauser, Vladimir Georgescu,  Thierry Jecko, Hubert Kalf,
Andreas Knauf, Michael Levitin, Fran\c cois Nicoleau, Heinz Siedentop
and Xue Ping Wang for  useful discussions. The first author was
partially supported by ESPRC grant EP/D054621. 

\section{Reduction of the problem}\label{s:red}
In this section, we study the resolvent of the
perturbed Dirac operator
\begin{eqnarray}\label{e:Hbd}
 H^{\rm bd}_\gamma=D_m+\gamma V, \mbox{ where } V:=v \otimes \Id_{\C^{2\nu}} 
 \mbox{ and $v$ \underline{bounded}}. 
\end{eqnarray} 
In Section \ref{s:mainresult}, we explain how to cover some
singularities. Due to the method, we will consider only small coupling
constants. We will show the limiting absorption principle 
\begin{eqnarray}\label{e:mainlow}
\sup_{|\lambda|\in [m, m+\delta], \, \varepsilon>0, \, |\gamma|\leq \kappa} \|\langle
Q\rangle^{-1}(H_\gamma^{\rm bd} -\lambda - \rmi \varepsilon )^{-1} \langle
Q\rangle^{-1}\|\leq C,
\end{eqnarray} 
for some $\kappa>0$. We notice this is equivalent to
\begin{eqnarray}\label{e:mainlow2}
\sup_{\lambda\in [m, m+\delta], \, \varepsilon>0, \, |\gamma|\leq \kappa} \|\langle
Q\rangle^{-1}(H_\gamma^{\rm bd} -\lambda - \rmi \varepsilon )^{-1} \langle
Q\rangle^{-1}\|\leq C,
\end{eqnarray} 
Indeed, by setting $\alpha_5:=\alpha_1\alpha_2\alpha_3\alpha_4$ and using the anti-commutation relation \eqref{e:Diracrep0}, we infer
\begin{eqnarray*}
\alpha_5\left(D_m+\gamma V\right)\alpha_5^{-1}=
-D_m+\gamma V.
\end{eqnarray*}
Note that $\alpha_5$ is unitary and stabilizes $\Hr^1(\R^3;
\C^{2n})$. This gives
\begin{eqnarray}\label{e:gauge}
\alpha_5 \varphi( D_m +\gamma V) \alpha_5^{-1}= 
 \varphi\big( -(D_m -\gamma V)\big), \mbox{ for all } \varphi\in \Cc(\R; \C).
\end{eqnarray}
  
\subsection{The non self-adjoint operator}\label{s:reduction}
Here, we relate the resolvent of \eqref{e:Hbd} in a point $z\in
\C\setminus \R$ with the one of some non self-adjoint Laplacian type  
operator $\Delta_{m,v,z}$, chosen in \eqref{e:Delta}. 
We fix a \emph{compact} set $\Ic$ being the area of energy we are 
concentrating on. In the next section, we explain how to recover a
limiting absorption principle for $H^{\rm bd}_\gamma$ over $\Ic$ given
the one of $\Delta_{m,\gamma v,z}$. 

We consider a potential $v\in L^\infty(\R^3;\R)$, not necessarily smooth,  satisfying  
\begin{equation}\label{a:BoundV}
\|v\|_\infty \leq m/2\mbox{ and }\nabla v\in
L^\infty(\R^3;\R^3).
\end{equation}
It particular, $(v(Q)-m-z)^{-1}$ stabilizes
$\Hr^1(\R^3; \C^{2\nu})$ for all $z$ in $\C\setminus \R$.

Since $\beta=\alpha_4$ satisfies \eqref{e:Diracrep0}, 
we deduct that $\beta$ has the eigenvalues $\pm 1$ and the eigenspaces 
have the same dimension. Let $P^+$ be the orthogonal
projection on the spin-up part of the space, i.e., on $\ker(\beta
-1)$. Let $P^-:= 1-P^+$. Since $\alpha_j$ satisfies \eqref{e:Diracrep0}, for $j\in \{1,2, 3\}$, we get $P^{\pm}\alpha_jP^{\pm}=0$. 
We set:
\begin{eqnarray*}
\alpha^+_j:= P^+ \alpha_j P^- \mbox{ and }  \alpha^-_j:= P^- \alpha_j
P^+, \mbox{ for } j\in\{1, 2, 3\}.
\end{eqnarray*} 
They are partial isometries:
\begin{eqnarray*}
\big(\alpha^+_j\big)^*= \alpha^-_j, \quad \alpha^+_j \alpha^-_j = P^+ 
\mbox{ and } \alpha^-_j \alpha^+_j = P^-, \mbox{ for } j\in\{1,2, 3\}.
\end{eqnarray*} 
The relation of anti-commutation \eqref{e:Diracrep0} gives:
\begin{eqnarray}\label{e:Diracrep1}
\alpha_i^-\alpha_j^+ + \alpha_j^-\alpha_i^+= 2\delta_{ij} \mbox{ and }
\alpha_i^+\alpha_j^- + \alpha_j^+\alpha_i^-= 2\delta_{ij}, \mbox{ for } i,j\in\{1,2, 3\}.
\end{eqnarray} 
We set $\C^\nu_\pm:= P^\pm\C^{2\nu}$.
In the direct sum $\C^{\nu}_-\oplus \C^{\nu}_+$, with a slight
abuse of notation, one can write 
\begin{eqnarray*}
\beta = \left(\begin{array}{cc}
\Id_{\C^\nu} & 0
\\
0 & -\Id_{\C^\nu}
\end{array}\right) \mbox{ and }
\alpha_j= \left(\begin{array}{cc}
0 & \alpha_j^+
\\
\alpha_j^- & 0
\end{array}\right),  \mbox{ for } j\in\{1,2, 3\}.
\end{eqnarray*}

We now split the Hilbert space $\Hr=L^2(\R^3; \C^{2\nu})$ into the spin-up
and down part: 
\begin{eqnarray}\label{e:split}
 \Hr=\Hr^+\oplus \Hr^-, \mbox{ where }  \Hr^\pm :=  L^2(\R^3;
 \C^{\nu}_\pm)\simeq L^2(\R^3;  \C^{\nu}) . 
\end{eqnarray} 
We define the operator:
\begin{equation}\label{e:Delta}
\Delta_{m,v,z}:=\alpha^+\!\cdot P\frac{1}{m -v(Q) +z}\alpha^-\!\cdot P 
+ v(Q)
\end{equation}
on $\Cc^\infty_c(\R^3; \C^{\nu}_+)$. It is well defined by
  (\ref{a:BoundV}). It is closable as its adjoint has a dense domain.
We consider the minimal extension, its closure. We denote its domain
by $\Dc_{\rm   min}(\Delta_{m,v,z})$ and keep the same symbol for the operator.
It is well known that even for symmetric operators one needs to be
careful with domains as the domain of the adjoint could be much bigger
than the one of the closure. In the next Proposition, we care about
this problem in our non-symmetric setting.

\begin{proposition}\label{p:domain}
Let $z\in \C\setminus \R$ such that $\Re(z)\geq 0$. Under the
hypotheses \eqref{a:BoundV}, we have that  
\begin{eqnarray*}
\Dc(\Delta_{m,v,z})= \Dc(\Delta_{m,v,\overline{z}}^* )= \Hr^2(\R^3;
\C^\nu_+) \mbox{ and } \Delta_{m,v,z}=\Delta_{m,v,\overline{z}}^*.
\end{eqnarray*} 
\end{proposition} 
\proof We mimic the Kato-Rellich approach and compare $\Delta_{m,v,z}$
with the more convenient operator $\tilde \Delta_z:=\big(1/(m+z)\big)
\Delta_{1,0,0}$. Its domain is $\Hr^2(\R^3; \C^\nu_+)$ and its
spectrum is $\{(m+\overline{z})t \mid t\in [0,\infty)\}$. We now show 
  there is $a\in [0,1)$ and $b\geq 0$ such that
\begin{eqnarray}\label{e:KR}
\|B f\|^2\leq  a\left\|\tilde\Delta_z f\right\|^2 +b\|f\|^2, 
\end{eqnarray} 
holds true for all $f\in \Cc^\infty_c(\R^3, \C^\nu_+)$,
where 
\begin{eqnarray*}
B:= \frac{v}{(m-v+z)} \tilde \Delta_z -i \frac{(\alpha^+ \cdot 
\nabla v)(Q)}{(m-v(Q)+z)^2}\,  \alpha^- \cdot P+ v(Q).
\end{eqnarray*}
Since $\|v\|_\infty \leq m/2$, $\Re(z)\geq 0$ and $\Im(z)>0$, we infer
$a_0:=\|v/(m-v+z)\|_\infty<1$. Set $M:=\|(\alpha^+\cdot  
\nabla v)(\cdot)/(m-v(\cdot)+z)^2\|_\infty$. Take $\varepsilon,
\varepsilon'\in (0,1)$
\begin{align*}
\|Bf\|^2 &\leq (1+\varepsilon)a_0^2\left\|\tilde \Delta_z f\right\|^2 + \left(1+ \frac{1}{\varepsilon}\right)
\left\|\frac{(\alpha^+ \cdot 
\nabla v)(Q)}{(m-v(Q)+z)^2}\,  \alpha^- \cdot P  f + v(Q )f\right\|^2, 
\\
&\leq (1+ \varepsilon) a_0^2\left\|\tilde \Delta_z f\right\|^2
+\frac{4M^2}{\varepsilon} \left\|\alpha^- \cdot P f\right\|^2 
+ \frac{4 \|v\|_\infty}{ \varepsilon^2 }\|f\|^2,
\\
&\leq \left((1+ \varepsilon) a_0^2 + \varepsilon'\right)
\left\|\tilde \Delta_z f\right\|^2 + \left(
\frac{4 \|v\|_\infty}{ \varepsilon^2 }+ \frac{2|m+z|^2M^2}{\varepsilon
 \varepsilon' } \right) \left\|f\right\|^2.
\end{align*} 
By choosing $\varepsilon$ and $\varepsilon' $ so that the first
constant is smaller than $1$, \eqref{e:KR} is
fulfilled.

Now, observe that $\|B(\tilde
\Delta_z+\mu)^{-1}\|^2\leq a +b\mu^{-2}$ for $\mu>0$. Fix $\mu_0>0$ such that 
$\|B(\tilde \Delta_z+\mu_0)^{-1}\|<1$. Then $(1+B(\tilde
\Delta_z+\mu_0)^{-1})$ is bijective. Noticing that
\begin{eqnarray*}
\big(\Id + B (\tilde \Delta_z+\mu_0)^{-1}\big)(\tilde \Delta_z+\mu_0) =
  \Delta_{m,v,z}+ \mu_0, 
\end{eqnarray*} 
we infer that $\Delta_{m,v,z}+ \mu_0$ is bijective from $\Hr^2(\R^3;
\C^\nu_+)$ onto $L^2(\R^3; \C^\nu_+)$. In particular
$\Dc(\Delta_{m,v,z})= \Hr^2(\R^3;\C^\nu_+)$. Directly, one has
$\Dc(\Delta_{m,v,z})\subset \Dc(\Delta_{m,v,\overline{z}}^*)$ and
$\Delta_{m,v,z}\subset \Delta^*_{m,v,\overline{z}}$ (inclusion of
graphs). Take now $f\in \Dc(\Delta_{m,v,\overline{z}}^* )$. Since
$\Delta_{m,v,z}+ 
\mu_0$ is surjective, there is $g\in \Dc(\Delta_{m,v,z})$ so that
\begin{eqnarray*}
(\Delta_{m,v,z}+ \mu_0) g = (\Delta_{m,v,\overline{z}}^* + \mu_0) f.
\end{eqnarray*} 
In particular, $(\Delta_{m,v,\overline{z}}^* + \mu_0)
(f-g)=0$. As   $\Delta_{m,v,z}+ 
\mu_0$ is surjective, we derive that $\ker
(\Delta_{m,v,\overline{z}}^*+\mu_0)=\{0\}$. 
In particular $f=g$, $\Dc(\Delta_{m,v,z})=
\Dc(\Delta_{m,v,\overline{z}}^*)$ and
$\Delta_{m,v,z}=\Delta^*_{m,v,\overline{z}}$. \qed

As a corollary, we derive:
\begin{lemma} The spectrum of $\Delta_{m,v,z}$ is
contained in the lower/upper half-plane which does not contain $z$. In
particular, $c+z$ is always in the resolvent set of $\Delta_{m,v,z}$
for any $c\in\R$. 
\end{lemma}
\proof Take now $f\in \Hr^2(\R^3; \C^\nu_+)$. Since
\begin{eqnarray}\label{e:imdelta}
\Im  \langle f, \Delta_{m,v,z}\,  f\rangle = \langle \alpha^-\!\cdot P f, 
\frac{-\Im(z)}{\big(m-v(Q)+\Re(z)\big)^2+\Im(z)^2}\, \alpha^-\!\cdot P
f\rangle, 
\end{eqnarray} 
is of the sign of $-\Im(z)$.  Since  $\Delta_{m,v,z}$ is a closed operator
having the same domain of its adjoint, the spectrum of
$\Delta_{m,v,z}$ is contained in the closure of its numerical range,
see Lemma \ref{l:NRT}.  \qed

We give a kind of Schur's lemma, so as to compute the inverse of the
Dirac operator. 
\begin{lemma}\label{l:inv} Suppose (\ref{a:BoundV}). Take $z\in
  \C\setminus \R$ such that $\Re(z)\geq 0$. In the 
spin-up/down decomposition $\Hr=\Hr^+\oplus \Hr^-$, 
we have $(H^{\rm bd}_1-z)^{-1}=$ 
\begin{align*}
\left(\begin{array}{c}
(\Delta_{m,v,z}+m-z)^{-1}
\\
\displaystyle
\frac{1}{m -v(Q) +z}\alpha^-\!\cdot P (\Delta_{m,v,z}+m-z)^{-1}
\end{array}\right. 
&
\\
&\hspace*{-4cm}
\left.\begin{array}{c}
\displaystyle
(\Delta_{m,v,z}+m-z)^{-1}\alpha^+\!\cdot P\frac{1}{m -v(Q) +z}
\\
\displaystyle
\frac{1}{m -v(Q) +z}\alpha^-\!\cdot P (\Delta_{m,v,z}+m-z)^{-1}
\alpha^+\!\cdot P\frac{1}{m -v(Q) +z}
-\frac{1}{m -v(Q) +z}
\end{array}\right) 
\end{align*}
\end{lemma}
\begin{remark}
The operator $(H^{\rm bd}_1-z)^{-1}$ is bounded from
$L^2(\R^3;\C^{2\nu})$ into $\Hr^1(\R^3;\C^{2\nu})$. However, this
improvement in the Sobolev scale does not hold if one looks at the
matricial terms separately. There is a real compensation coming from the
off-diagonal terms. First note that $\alpha^- \cdot
P(\Delta_{m,v,z}+m-z)^{-1} \alpha^+\cdot P$ is a bounded 
operator in $L^2(\R^3; \C^\nu_-)$ and a priori not into 
$\Hr^s(\R^3; \C^\nu_-)$, with $s>0$. Indeed, $\alpha^+\cdot P$ sends
$L^2(\R^3; \C^\nu_-)$ into $\Hr^{-1}(\R^3; \C^\nu_+) $, then
$(\Delta_{m,v,z}+m-z)^{-1}$ to $\Hr^1(\R^3; \C^\nu_+)$ and the left
$\alpha^-\cdot P$ sends again into $L^2(\R^3; \C^\nu_-)$. 
On the other hand, the term $(\Delta_{m,v,z}+m-z)^{-1}$ is bounded
from $L^2(\R^3; \C^\nu_+)$ into $\Hr^2(\R^3; \C^\nu_+)$, which is much
better than expected.
\end{remark}
 
\begin{proof} 
Let $f\in L^2(\R^3;\C^{2\nu})$. By self-adjointness 
of $H^{\rm bd}_1$, there is a unique $\psi\in\Hr^1(\R^3;\C^{2\nu})$ such  that
$(H^{\rm bd}_1-z)\psi=f$. We separate the upper and lower spin components and
denote  $f=(f_+, f_-)$ and $\psi=(\psi_+, \psi_-)$ in $\Hr=\Hr^+\oplus
\Hr^-$.  We rewrite the equation $(D_m +V(Q)-z)\psi =f$ to get: 
\begin{equation}\label{e:system}
\begin{cases}
\alpha^+\!\cdot P \psi_- + m \psi_+  + v(Q)\psi_+ -z\psi_+=f_+,\\
\alpha^-\!\cdot P \psi_+ - m \psi_-  +v(Q)\psi_- -z\psi_-=f_-.
\end{cases}
\end{equation}
From the second line, we get $(v(Q)-m-z)\psi_-=f_- -\alpha^-\!\cdot P
\psi_+$. Since $z$ is not real, we can take the inverse and infer
$\psi_-=(v(Q)-m-z)^{-1}(f_- -\alpha^-\!\cdot P \psi_+)$. Since $\psi_- \in
\Hr^{1}$, we can apply it $\alpha^+\!\cdot P$ and obtain a vector
 of $L^2(\R^3;\C^\nu_+)$. Now, since $f_-$ is  in
 $L^2(\R^3;\C^\nu_-)$ and since  $(v(Q)-m-z)^{-1}$ is bounded, we have
 $\alpha^+\!\cdot P  (v(Q)-m-z)^{-1}f_-\in  \Hr^{-1}(\R^3;\C^\nu_+)$ and since
 $(v(Q)-m-z)^{-1}\alpha^-\cdot P  
\psi_+$ is in $L^2(\R^3;\C^\nu_-)$, we rewrite the system: 
\begin{equation*}
\begin{cases}
  \displaystyle \left(\alpha^+\!\cdot P\frac{1}{m -v(Q)
    +z}\,\alpha^-\!\cdot P   +v(Q) + m-z\right)\psi_+= &\displaystyle
  f_+ +\alpha^+\!\cdot P\frac{1}{m 
  -v(Q) +z}f_-,\\ 
\hfill\psi_-=&\displaystyle \frac{1}{m -v(Q) +z}\left(\alpha^-\!\cdot P
\psi_+ -f_-. \right)  
\end{cases}
\end{equation*}
To conclude it remains to show that $\Delta_{m,v,z}+m-z$ is invertible in
$\Bc(\Hr^{1}, \Hr^{-1})$, so as to invert it in the system. 
Using \eqref{e:imdelta}, we have $|\Im\langle u,
(\Delta_{m,v,z}-z) u \rangle| \geq c \|u\|_{\Hr^1}^2$. Then
$\|(\Delta_{m,v,z}+m-z) u\|_{\Hr^{-1}}\geq  c\|u\|_{\Hr^1} $ and 
$\|(\Delta_{m,v,z}+m-z)^* u\|_{\Hr^{-1}}\geq  c\|u\|_{\Hr^1}$ hold. Thus,
$\Delta_{m,v,z}-z$ is bijective from $\Hr^1$ 
onto $\Hr^{-1}$.
\end{proof} 
\subsection{From one limiting
absorption principle to another}\label{s:LAPano}
The main motivation for the operator $\Delta_{m,v,z}$ is to deduce a limiting
absorption principle for $H^{\rm bd}_\gamma$ starting with one for
$\Delta_{m,\gamma v,z}$. Consider the upper right term in Lemma
\ref{l:inv}, the basic idea would be to put by force the weight 
$\langle Q\rangle^{-1}$ and to say that every terms are bounded. However,
we have that 
\begin{eqnarray*}
\underbrace{\langle Q\rangle^{-1} (\Delta_{m,v,z}+m-z)^{-1} \langle Q
  \rangle^{-1}}_{{\rm bounded\, from\, LAP\, for\, }  \Delta_{m,v,z}}\,
\underbrace{\langle Q\rangle  \alpha^+\cdot P \langle
  Q\rangle^{-1}}_{\rm unbounded}\, \frac{1}{m -v(Q) +z}. 
\end{eqnarray*} 
One needs to take advantage that one seeks an estimate on a
bounded interval of the spectrum. Therefore, we start with a lemma of
localization in the momentum space and elicit a solution in Lemma \ref{l:LAPchange}. 
Note also that one may consider $\Im z <0$ by taking the adjoints
in the two next lemmata. We shall also use estimates which are uniform
in the coupling constant, due to Proposition \ref{p:collage}.

\begin{lemma}\label{l:steptoLAP}
Set $\Ic\subset \R$ a compact interval. Let $V$ be a bounded potential
and $\kappa>0$.  There is an \emph{even} function
$\varphi\in\Cc^\infty_c(\R;\R)$ such that   
the following estimations of the resolvent are equivalent:
\begin{equation}\label{e:LAPchange1}
\sup_{\Re z\in \Ic, \Im z>0, |\gamma|\leq \kappa}
\left\|\langle Q\rangle^{-1} \varphi(\alpha\cdot P) (D_m+\gamma V(Q)-z)^{-1}
\varphi(\alpha\cdot P) 
\langle Q\rangle^{-1}\right\|<\infty,  
\end{equation}
\begin{equation}\label{e:LAPchange1''}
\sup_{\Re z\in \Ic, \Im z>0, |\gamma|\leq \kappa}
\left\|\langle Q\rangle^{-1}(D_m+\gamma
V(Q)-z)^{-1}\varphi(\alpha\cdot P)\langle
Q\rangle^{-1}\right\|<\infty, 
\end{equation}
\begin{equation}\label{e:LAPchange1'}
\sup_{\Re z\in \Ic, \Im z>0, |\gamma|\leq \kappa}
\left\|\langle Q\rangle^{-1}(D_m+\gamma V(Q)-z)^{-1}\langle
Q\rangle^{-1}\right\|<\infty. 
\end{equation}
\end{lemma} 
\proof It is enough to consider $\Im z \in (0,1]$. Set
$\Jc:=\Ic\times(0,1]\times [-\kappa,\kappa]$, $H_\circ:=\alpha\cdot
P$ and $H_\gamma :=D_m+\gamma V$. We choose $\varphi_1\in\Cc^\infty_c(\R)$ with 
value in $[0,1]$, being even and equal  to $1$ in a neighborhood of
$0$. We define $\varphi_R(\cdot):=\varphi_1(\cdot/R)$ and $\widetilde \varphi_R:=1- \varphi_R$. 

We first notice that $\langle Q\rangle \in \Cc^1(H_\circ)$, see Appendix
\ref{s:dev-commut}. There is a constant $C>0$ so that 
\begin{eqnarray}\label{e:kommu}
\left|\langle \langle Q\rangle f, \alpha\cdot P f \rangle - \langle
\alpha \cdot P f, \langle Q\rangle 
f\rangle\right| = \left|\langle f, (\alpha\cdot \nabla \langle \cdot
\rangle)(Q) f\rangle\right|\leq  C\|f\|^2     
\end{eqnarray} 
holds true for all $f \in \Cc^\infty_c (\R^3 ;\C^{2\nu})$. This is usually not
enough to deduce the $\Cc^1$ property, see \cite{GeorgescuGerard}. We use
\cite{GoleniaMoroianu}[Lemma A.2] with the notations $A:=H_\circ$,
$H:=\langle Q\rangle$, $\cchi_n(x):=\varphi(x/n)$ and with
$\Dr:=\Cc^\infty_c (\R^3;\C^{2\nu})$. The hypotheses are fulfilled and we
deduce that  $\langle Q\rangle \in \Cc^1(H_\circ)$. 

By the resolvent equality,  we have:  
\begin{align}
\nonumber
(H_\gamma-z)^{-1}\tilde \varphi_R(H_\circ)\big(\Id + W( H_\circ
-z)^{-1}\tilde \varphi_R(H_\circ)\big)=&\, ( H_\circ-z)^{-1}
\tilde \varphi_R(H_\circ) 
\\ \label{e:LAPchange3}
&-(H_\gamma-z)^{-1} \varphi_R( H_\circ) W( H_\circ-z)^{-1} \tilde
\varphi_R( H_\circ),  
\end{align} 
where $W:=\gamma V+m\beta$.  Note that  the support of $\tilde \varphi_R$
vanishes as $R$  goes to infinity. We have
\begin{eqnarray*}
 \|\langle Q\rangle ( H_\circ -z)^{-1}\tilde \varphi_R( H_\circ) \langle Q
 \rangle^{-1}\|\leq \Oc(1/R), \mbox{  uniformly in } (z, \gamma)\in \Jc.
\end{eqnarray*}
 Indeed, if we commute with $\langle Q\rangle$, the part in $( H_\circ
 -z)^{-1}\tilde  \varphi_R( H_\circ)$ is a $\Oc(1/R)$ by functional
 calculus. For the   other part, Lemma  \ref{l:est3} gives 
\begin{eqnarray*}
\|[\langle Q\rangle, ( H_\circ -z)^{-1}\tilde \varphi_R( H_\circ)] \langle
Q\rangle^{-1}\|\leq \Oc(1/R^2), \mbox{  uniformly in } (z, \gamma)\in \Jc.
\end{eqnarray*}
Remembering $V$ is bounded and choosing $R$ big enough, we infer there
is a constant $c\in (0,1)$, so that 
\begin{eqnarray}\label{e:LAPchangeagain}
 \|W 
\langle Q\rangle(H_\circ-z)^{-1}\tilde\varphi_R( H_\circ)\langle
Q\rangle^{-1}\|\leq c, \mbox{  uniformly in } (z, \gamma)\in \Jc.
\end{eqnarray} 
We fix $R$ and choose $\varphi:=\varphi_R$. 
We now prove the equivalence. Observe that $\langle
Q\rangle^{-1}\varphi(H_\circ)\langle Q\rangle$ is bounded, 
since $\langle Q\rangle \in \Cc^1(H_\circ)$. 
One infers directly that \eqref{e:LAPchange1'} $\Rightarrow$
\eqref{e:LAPchange1''} $\Rightarrow$ \eqref{e:LAPchange1}. 
It remains to prove \eqref{e:LAPchange1} $\Rightarrow$
\eqref{e:LAPchange1'}. Thanks to \eqref{e:LAPchangeagain}, we deduce
from \eqref{e:LAPchange3} that:
\begin{align}\nonumber
\langle Q\rangle^{-1}(H_\gamma-z)^{-1}\tilde \varphi(H_\circ)\langle
Q\rangle^{-1}=&\, 
\Big(\langle Q\rangle^{-1}( H_\circ-z)^{-1} \tilde \varphi(
H_\circ)\langle Q\rangle^{-1} 
\\ \label{e:LAPchange2}
& -\langle Q\rangle^{-1}(H_\gamma-z)^{-1} \varphi( H_\circ) \langle
Q\rangle^{-1} \quad  W
\langle Q \rangle    ( H_\circ-z)^{-1} \tilde 
\varphi( H_\circ)\langle Q\rangle^{-1}\Big)
\\ \nonumber
&\times  \big(\Id + W\langle Q\rangle(H_\circ -z)^{-1}\tilde
\varphi(H_\circ) \langle Q\rangle^{-1}\big)^{-1}.   
\end{align}
Note that the last line and the right part of the second line of the
r.h.s.\ are uniformly bounded in $(z, \gamma)\in \Jc$ by
\eqref{e:LAPchangeagain}. We multiply on the left by the bounded
operator $\langle Q\rangle^{-1}\varphi(H_\circ)\langle Q\rangle$.  
The first term of the r.h.s.\ is bounded uniformly by functional
calculus. For the second one, we use \eqref{e:LAPchange1}. We infer:
\begin{equation*}
\sup_{(z, \gamma)\in \Jc}
\left\|\langle Q\rangle^{-1} \varphi (H_\circ)   (H_\gamma-z)^{-1}
\tilde  \varphi (H_\circ)
\langle Q\rangle^{-1}\right\|<\infty.
\end{equation*}
Doing like in \eqref{e:LAPchange2}, on the left hand side, we get
\begin{equation} \label{e:LAPchange4}
\sup_{(z, \gamma)\in \Jc}
\left\|\langle Q\rangle^{-1}  \tilde \varphi(H_\circ)
(H_\gamma-z)^{-1} \varphi(H_\circ)
\langle Q\rangle^{-1}\right\|<\infty.  
\end{equation}
Finally, to control $\langle Q\rangle^{-1}  \tilde\varphi (H_\circ)
(H_\gamma-z)^{-1} \tilde\varphi (H_\circ) \langle Q\rangle^{-1}$, we
multiply 
\eqref{e:LAPchange2} on the left by the bounded operator
$\langle Q\rangle^{-1}\tilde\varphi (H_\circ) \langle Q\rangle$ and deduce 
the boundedness using \eqref{e:LAPchange4}.\qed 

\begin{lemma}\label{l:LAPchange}
Take $\kappa\in (0,1]$ and  a compact interval $\Ic \subset [0,\infty)$. Suppose (\ref{a:BoundV}) and that
\begin{equation}\label{e:LAPchange0}
\sup_{\Re z\in \Ic,\Im z\in (0,1], |\gamma|\leq \kappa}
\left\|\langle Q\rangle^{-1}(\Delta_{m,\gamma v,z}+m-z)^{-1}\langle
Q\rangle^{-1}\right\|<\infty 
\end{equation}
hold true. Then, we have
\begin{equation}\label{e:LAPchange}
\sup_{\Re z\in \Ic, \Im z>0, |\gamma|\leq \kappa}
\left\|\langle Q\rangle^{-1}(D_m+ \gamma v(Q)\otimes 
\Id_{\C^{2\nu}} -z)^{-1}\langle Q\rangle^{-1}\right\|<\infty.
\end{equation}
\end{lemma}
\proof Set $H_\circ:=\alpha\cdot P$ and $ \Jc:=\Ic\times
  (0,1]\times [-\kappa, \kappa]$. 
By Lemma \ref{l:steptoLAP}, it is enough to show 
\eqref{e:LAPchange1} for a chosen $\varphi$. Since $\varphi$ is even
and constant in a neighborhood of $0$, by setting
$\psi(\cdot):=\varphi(\sqrt{|\cdot|})$, we have $\psi\in\Cc^\infty_c(\R)$
and that $\varphi(H_\circ)= \psi\big((\alpha\cdot P)^2 \big)$. 
In particular, we obtain $\varphi(H_\circ)$ stabilizes $\Hr^{\pm}$ and
have the right to let it appear in spin decomposition of the resolvent
of $H$ of Lemma \ref{l:inv}. We treat only the upper right corner of the
expression as the others are managed in the same way. We need to
bound the term:  
\begin{align*}
\langle Q\rangle^{-1} \varphi(H_\circ)
(\Delta_{m,\gamma v,z}+m-z)^{-1}\alpha^+ \cdot P\frac{1}{m 
-\gamma v(Q) +z} \varphi(H_\circ)\langle Q\rangle^{-1} &=
\\
&\hspace*{-9cm}\langle Q\rangle^{-1} \varphi(H_\circ)\langle Q\rangle\quad 
\langle Q\rangle^{-1}(\Delta_{m,\gamma v,z}+m-z)^{-1}\langle Q\rangle^{-1}
\quad \langle Q\rangle \alpha^+ \cdot P\frac{1}{m  
-\gamma v(Q) +z} \varphi(H_\circ)\langle Q\rangle^{-1}. 
\end{align*} 
The middle term is controlled by the hypothesis. Thanks to
\eqref{e:kommu}, one has that $[\varphi(H_\circ),\langle Q\rangle]$ is bounded;
the first term is bounded. For the last one, we commute:
\begin{align*}
 \langle Q\rangle \alpha^+ \cdot P\frac{1}{m -\gamma v(Q) +z}
 \varphi(H_\circ)\langle Q\rangle^{-1} =&  
\langle Q\rangle \left[\alpha^+ \cdot P,\frac{1}{m -\gamma v(Q)
    +z}\right] \langle Q\rangle^{-1}\, \quad \langle Q\rangle 
\varphi(H_\circ)\langle Q\rangle^{-1} 
\\ &+ 
\langle Q\rangle\frac{1}{m -\gamma v(Q) +z} \langle Q\rangle^{-1}\, \quad
\langle Q\rangle\alpha^+ \cdot P\varphi(H_\circ)\langle Q\rangle^{-1}.  
\end{align*} 
We estimate uniformly in $(z, \gamma)\in \Jc$.
By \eqref{a:BoundV} , we get $\|\langle Q\rangle (m-\gamma v(Q)+z)^{-1} \langle
Q\rangle^{-1}\|$ is bounded as $\langle Q\rangle$ commute 
with $v$. By \eqref{a:BoundV}, we also obtain that $\|\langle Q\rangle
[\alpha^+ \cdot P,   (m-\gamma v+z)^{-1}]\langle Q\rangle^{-1}\|$ is also
controlled. At last, it is enough to consider $\langle
Q\rangle \partial_j \varphi(H_\circ)\langle Q\rangle^{-1}$, which is
easily bounded by Lemma \ref{l:est3} for instance. \qed

We come to other types of weights. Motivated by the non-relativistic
case, see Theorem \ref{t:nonrelamain}, we are interested in  
singular weights like $|Q|$. But, as noticed in Remark \ref{r:noLAP},
the operator $|Q|^{-1}(H^{\rm bd}-z)^{-1} |Q|^{-1}$ is even not
bounded. Therefore, we enlarge the space in momentum and try the first
reasonable weight, namely $\langle P\rangle^{1/2}|Q|$. Given $z\in
\C\setminus \R$ and using the Hardy inequality, one reaches 
\begin{align*}
\|\langle P\rangle^{-1}|Q|^{-1} ( H^{\rm bd}_\gamma -v -z)^{-1}|Q|^{-1}\|&\leq 
\| \langle P\rangle^{-1} |P|\, \| \cdot \|\, |P|^{-1} |Q|^{-1}\|^2 \cdot
\|( H^{\rm bd}_\gamma -v -z)^{-1} |P|\,\| 
\\
&\leq \, C(\kappa) \langle z \rangle / |\Im(z)|.
\end{align*}
By interpolation, one infers 
\begin{align*}
\|\langle P\rangle^{-1/2}|Q|^{-1} ( H^{\rm bd}_\gamma -v -z)^{-1}|Q|^{-1}\langle P\rangle^{-1/2}\| \leq C(\kappa) \langle z \rangle / |\Im(z)|<\infty.
\end{align*}
The upper bound seems relatively sharp in $z$. However, under the same hypotheses as before, we obtain:
\begin{lemma}\label{l:LAPchange2}
Take $\kappa\in (0,1]$ and  a compact interval $\Ic \subset [0,\infty)$. Suppose (\ref{a:BoundV}) and that \eqref{e:LAPchange0} hold true. Then, there is $C>0$ so that
\begin{equation}\label{e:LAPchange'}
\sup_{\Re z\in \Ic, \Im z >0, |\gamma|\leq \kappa}
|\langle f, (D_m+ \gamma v(Q)\otimes 
\Id_{\C^{2\nu}} -z)^{-1} f \rangle|\leq C \|  \langle P\rangle^{1/2}|Q| f\|^2.
\end{equation}
\end{lemma}
\proof It is enough to consider $\Im z \in (0,1]$.
Set $H_\gamma :=D_m+\gamma v(Q)\otimes  \Id_{\C^{2\nu}}$ and
$\Jc:=\Ic\times(0,1]\times [-\kappa,\kappa]$. 
Let $f=(f_+,f_-)$, with $f_{\pm}\in \Cc^\infty_c(\R^3\setminus\{0\}; \C^{\nu}_{\pm})$. Lemma \ref{l:inv} and \eqref{a:BoundV} give a constant $C>0$, uniform in $(z, \gamma)\in \Jc$, so that:
\begin{align*}
\left|\left\langle f, (H_\gamma -z)^{-1}f\right\rangle\right|\leq&
\,4/m^2 \|f_-\|^2 + 2 \left\||Q|^{-1}(\Delta_{m,\gamma v,z}+m-z)^{-1}|
Q|^{-1}\right\| \times
\\
&\hspace*{-2cm}\left(\big\||Q|f_+\big\|^2+\big\||Q|\alpha^+\cdot P (m-v(Q)+z)^{-1}f_-\big\|^2+\big\||Q|\alpha^+\cdot P(m-v(Q)+\overline{z})^{-1}f_-\big\|^2\right) 
\\
\leq& \,C \left(\big\||Q|f_+\big\|^2 + \big\||Q|\alpha^+\cdot P
f_-\big\|^2 +  \|f_-\|^2 \right).  
\end{align*}
Note that the Hardy inequality gives that $\|f_-\|\leq 2 \big\| |Q|
\alpha^+\cdot P f_- \big\|$. Then, by commuting $|Q|$ with
$\alpha^+\cdot P$ over $\Cc^\infty_c(\R^3\setminus\{0\}; \C^\nu)$, we
find $C'>0$, so that: 
\begin{align*}\sup_{(z, \gamma)\in \Jc }\left|\left\langle  f,
  |Q|^{-1}(H_\gamma -z)^{-1} |Q|^{-1}f\right\rangle\right|\leq& C'
  \left(\big\|\Id\otimes P_+\, f\big\|^2 + \big\|\langle P\rangle
  \otimes P_-\, f\big\|^2\right), 
	\end{align*}
for all $f\in \Cc^\infty_c(\R^3\setminus\{0\}; \C^{2\nu})$.
Here we identify, $L^2(\R^3; \C^{2\nu})\simeq L^2(\R^3)\otimes \C^{2\nu}$. 
We now exchange the role of $P_+$ and $P_-$. Considering the operator  
\[\alpha^-\!\cdot P\frac{1}{m +v(Q) +z}\alpha^+\!\cdot P -v(Q) \mbox{
  in } L^2(\R^3; \C^{\nu}_-),\] which leads to the same arguments as
for $\Delta_{m,-v,z}$ if one identifies $\C^{\nu}_-\simeq \C^{\nu}_+$,
one obtains also that  
\begin{align*}
\left|\left\langle f, |Q|^{-1}(H_\gamma -z)^{-1}|Q|^{-1}
f\right\rangle\right|\leq& C' \left(\big\|\Id\otimes P_-\, f\big\|^2 +
\big\|\langle P\rangle \otimes P_+\, f\big\|^2\right),  
\end{align*}
for all $f\in \Cc^\infty_c(\R^3\setminus\{0\}; \C^{2\nu})$. By
interpolation, e.g., \cite{BerghLofstrom}[Theorem 4.4.1 and Theorem
  6.4.5.(7)], we infer:  
\begin{align*}
\left|\left\langle  f, |Q|^{-1}(H_\gamma -z)^{-1}
|Q|^{-1}f\right\rangle\right|\leq& C'' \big\|\langle P\rangle^{1/2}
f\big\|^2, 
\end{align*}
for all $f\in \Cc^\infty_c(\R^3\setminus\{0\}; \C^{2\nu})$. The latter is a core 
in $\Hr^1(\R^3; \C^{2\nu})$.\qed

\section{Positive commutator estimates.}\label{s:positive}
In the previous section, we saw how to deduce some estimate of the
resolvent for $D_m+V(Q)$ starting with some of $\Delta_{m,v,z}$,
namely \eqref{e:LAPchange0}. First, one technical problem is that operator
depends on the spectral parameter; hence we will study a family of
operators uniformly in the spectral parameter. Secondly, we are
concerned about the interval $[m, m+\delta]$ and we know that there is
no such estimate above $(m- \varepsilon , m)$ as eigenvalues usually
accumulates to $m$ from below. 
Since the theory developed in Appendix \ref{s:WeakMourreTheory} gives
some estimates for $\Re z\in [0, \infty)$, we will perform a
  shift. Therefore, we study the operator 
\begin{eqnarray*}
\Delta_{2m,\gamma v,\xi}, \mbox{ uniformly in } (\gamma, \xi)\in \Ec= \Ec(\kappa, \delta):=
      [-\kappa,\kappa]\times [0,\delta]\times(0,1].
\end{eqnarray*}
Here we use a slight abuse of notation identifying $\C\simeq
\R^2$. One should read $\Re\xi\in [0, \delta], \Im\xi\in(0, 1]$ 
and $|\gamma|\leq \kappa$. Note the uniformity in the coupling
constant is used in Proposition \ref{p:collage}.

To show \eqref{e:LAPchange0}, and therefore \eqref{e:LAPchange1'} with
the help of Lemma \ref{l:steptoLAP}, it is enough to prove the
following fact. Note we strengthen the hypothesis \eqref{a:BoundV}. 
\begin{theorem}\label{t:mourreD}
Suppose that $v\in L^\infty(\R^3; \R)$ satisfies the hypothesis (H1)
and (H2) from Theorem \ref{t:vraimain}. Then there are $\delta, \kappa, C_{\rm
  LAP} >0$ such that  
\begin{eqnarray}\label{e:mourreD}
\sup_{\Re z\geq 0,\Im z>0, (\gamma, \xi)\in \Ec}
|\langle f, (\Delta_{2m,\gamma v,\xi}-z)^{-1} f\rangle | \leq C_{\rm
  LAP} \big\|\, |Q| f \big\|^2.
\end{eqnarray}
\end{theorem} 

We will show the theorem in the end of the section. We proceed by
checking the hypothesis of Appendix \ref{s:WeakMourreTheory}.  We
recall \eqref{e:Diracrep1} and fix some notation: 
\begin{eqnarray*}
S:= \Delta_{1,0,0}= \alpha^+ \cdot P\, \alpha^- \cdot P = -\Delta_{\R^3} \otimes
\Id_{\C^\nu_+} \mbox{ in } 
\Hr^2(\R^3; \C^\nu_+) \simeq \Hr^2(\R^3)\otimes \C^\nu_+
\end{eqnarray*} 
and set $\Sr:=\dot\Hr^1(\R^3; \C^\nu_+)$, 
the homogeneous Sobolev space of order $1$, i.e., the completion of 
$\Hr^1(\R^3; \C^\nu_+)$ under the norm $\|f\|_\Sr:= \| S^{1/2} f\|^2$. 
Consider the strongly continuous one-parameter unitary group
$\{W_t\}_{t\in \R}$ acting by:
\begin{eqnarray*}
(W_t f)(x)= e^{3t/2} f(e^tx), \mbox{ for all } f\in L^2(\R^3; \C^\nu_+).
\end{eqnarray*} 
This is the $C_0$-group of dilatation. Easily, by interpolation and
duality, one gets  
\begin{eqnarray*}
W_t \Sr \subset \Sr \mbox{ and } W_t \Hr^s(\R^3; \C^\nu_+) \subset
\Hr^s(\R^3; \C^\nu_+), \mbox{ for all } s\in \R.
\end{eqnarray*} 
Consider now its generator $A$ in $L^2(\R^3; \C^\nu_+)$. 
By the Nelson lemma, it is essentially self-adjoint on
$\Cc^\infty_c(\R^3; \C^\nu_+)$. It acts as follows:
\begin{eqnarray*}
\displaystyle A=\frac{1}{2}(P\cdot Q+ Q\cdot P)\otimes \Id_{\C^\nu_+} \mbox{ on }
\Cc^\infty_c(\R^3; \C^\nu_+) \simeq \Cc^\infty_c(\R^3)\otimes \C^\nu_+. 
\end{eqnarray*}
In the next Proposition, we will choose the upper bound $\kappa$ of
the coupling constant and state the commutator estimates. 

\begin{proposition}\label{p:commu1}
Let $\delta\in (0, 2m)$. Suppose that the hypotheses (H1) and (H2) are
fulfilled. Then there are $c_1, \kappa>0$ such that 
\begin{align}
\label{e:commu10}
\Dc(\Delta_{2m,\gamma v,\xi}) = \Hr^2(\R^3; \C^\nu_+), \quad
(\Delta_{2m,\gamma v,\xi})^* = \Delta_{2m,\gamma v, \overline{\xi}},
\\
\label{e:commu11}
[\Re(\Delta_{2m,\gamma v,\xi}),\rmi A]_\circ- c_v\Re(\Delta_{2m, \gamma v,\xi}) \geq
c_1 S>0,  
\\
\mp[\Im (\Delta_{2m, \gamma v, \Re(\xi) \pm\rmi \Im(\xi)}) , \rmi A]_\circ\geq 0,
\end{align} 
hold true in the sense of forms on $\Hr^1(\R^3; \C^\nu_+)$, for all
$(\gamma, \xi)\in \Ec$. 
\end{proposition}

\proof
We start with a first restriction on $\kappa$. We set $\kappa\leq
(2m-\delta) / \|v\|_\infty$. Hence, 
\begin{eqnarray}\label{e:bd}
\delta \leq 2m - \gamma v(\cdot) + \Re(\xi) \leq 4m, \mbox{ for all }
(\gamma, \xi)\in \Ec.
\end{eqnarray} 
In particular, $0$ is not in the essential image of $2m- \gamma v +\Re(\xi)$;
Proposition \ref{p:domain} gives \eqref{e:commu10}. 

We turn to the commutator estimates. It is enough to compute
the commutators in the sense of form on  $\Cc^\infty_c(\R^3; \C^\nu)$,
since it is a core for $\Delta_{2m,v,\xi}$ and $A$. 
\begin{align}
\nonumber
\big[\Delta_{2m,\gamma v,\xi},\rmi A\big]
=& \left[\alpha^+ \cdot P\, \frac{1}{2m -\gamma v +\xi}\alpha^- \cdot 
P,\rmi A\right] +\gamma [v,\rmi A]\\ 
\label{e:firstcommu}
=&\, 2\,\alpha^+\cdot P\, \frac{1}{2m -\gamma v +\xi}\alpha^- \cdot P
-\gamma\, \alpha^+\cdot P\, \frac{Q\cdot \nabla v(Q)}{(2m - \gamma v
+\xi)^2}\alpha^-\cdot P  - \gamma\, Q\cdot\nabla v(Q).
\end{align}
Then, we have $[\Re(\Delta_{2m,\gamma v,\xi}), \rmi A]-
c_v\Re(\Delta_{2m, \gamma v,\xi})=$ 
\begin{align*}
=&\,(2-c_v)\, \alpha^+\cdot P\frac{2m -\gamma v+\Re(\xi)}{\big(2m 
-\gamma v+\Re(\xi)\big)^2+\Im(\xi)^2}\alpha^-\cdot P
\\
& 
-\gamma\,\alpha^+\cdot P\left(\frac{Q\cdot \nabla  v(Q)\big(\big(2m 
-\gamma v+\Re(\xi)\big)^2-\Im(\xi)^2\big)}{\big(\big(2m - \gamma
  v+\Re(\xi)\big)^2+\Im(\xi)^2\big)^2 }\right) 
\alpha^-\cdot P
 - \gamma\, Q\cdot\nabla  v(Q) - c_v \gamma v(Q).
\\
\geq & (2-c_v)\frac{\delta}{16 m^2+1} S - \kappa\, \|Q\cdot \nabla v(Q)\|\,
\frac{ 16 m^2 +1}{ \delta^4} S -\kappa \frac{c_v'}{|Q|^2} \geq c_1 S,
\end{align*}
where $\displaystyle c_1:= \frac{\delta (2-c_v) }{32 m^2+2}$ and by
assuming that 
$\displaystyle \kappa \leq \frac{c_1}{2 (4 c_v'+ \|Q\cdot \nabla v(Q)\|
( 16 m^2 +1)/\delta^4)}$. Note the ``$4$'' comes from the Hardy
inequality.  This gives \eqref{e:commu11}. 
 
At last, we have:
\begin{align*}
[\Im\Delta_{2m, \gamma v,\xi}, \rmi A] &= -2 \Im(\xi)\, \alpha^+\cdot P
\frac{\big(2m-\gamma v+\Re(\xi)\big)^2+\Im(\xi)^2 - \gamma Q\cdot
  \nabla  v(Q)\big(2m-\gamma v+\Re(\xi)\big) }
{\big(\big(2m-\gamma v+\Re(\xi)\big)^2+\Im(\xi)^2\big)^2} 
\alpha^- \cdot P.
\end{align*}
This is of the sign of $-\Im(\xi)$, when we further impose
$\kappa \leq \delta^2/(8 m \|Q\cdot \nabla  v(Q)\|)$. \qed

We now bound some commutators. 
\begin{proposition}\label{p:commu2}
Let $\delta\in (0, 2m)$. Suppose that the hypotheses (H1) and (H2) are
fulfilled. Consider the $c_1, \kappa>0$ from Proposition
\ref{p:commu1}. There is $c$ and $C$ depending on $c_v, \delta, \kappa$ and
$v$, such that
\begin{eqnarray}\label{e:mourrestrict0chk}
|\langle \Delta_{2m,\gamma v,\overline{\xi}}\, f, Ag \rangle - \langle
Af, \Delta_{2m,\gamma v,\xi}\, g \rangle   |\leq  
c\|f\|\cdot\|(\Delta_{2m,\gamma v,\xi} \pm i) g\|, 
\end{eqnarray}  
holds true, for all $f,g \in \Hr^2(\R^3; \C^+_\nu)\cap \Dc(A)$
and 
\begin{align}
\label{e:commu21}
|\langle f, [[\Delta_{2m,\gamma v,\xi},\rmi A]_\circ, \rmi A]_\circ
f\rangle| \leq C  \langle f, S f\rangle. 
\end{align}
holds true for all $f\in\Hr^{1}(\R^3; \C^\nu_+)$.
\end{proposition}
\proof We take $\kappa$ as in the proof of Proposition
\ref{p:commu1}. We first find $c>0$, uniform in 
$(\gamma, \xi)\in \Ec$, so that
\begin{eqnarray}\label{e:commuS}
|\langle f,  [ \Delta_{2m, \gamma v, \xi }, iA ]_\circ g\rangle | \leq c \big(\| f
\|\cdot\|g\| + \|f\|\cdot \|Sg\|\big), \mbox{ for all }  f,g \in
\Cc^\infty_c(\R^3; \C^\nu_+).
\end{eqnarray} 
We recall that the latter space is a core for $A$, $S$ and $\Delta_{2m,
  \gamma v, \xi }$. Taking in account \eqref{e:firstcommu}, observe that
\begin{eqnarray*}
\left|\frac{2}{2m -\gamma v +\xi}
-\gamma \frac{Q\cdot \nabla v(Q)}{(2m - \gamma v
+\xi)^2}\right|\leq \frac{2}{\delta} + \kappa \frac{\|Q\cdot \nabla
  v(Q)\|}{\delta^2}. 
\end{eqnarray*} 
It remains to find $a,b>0$, which are uniform in $(\gamma, \xi)\in \Ec$, such
that the following estimation holds:
\begin{eqnarray*}
\|\Delta_{2m, \gamma v, \xi } f\|\geq a \|S f\| - b\|f\|, \mbox{ for
  all } f\in \Cc^\infty_c(\R^3, \C^\nu_+).
\end{eqnarray*} 
This follows from $\|\Delta_{2m, \gamma v, \xi } f\|^2\geq a^2 \|S
f\|^2 - b^2\|f\|^2$. Take
$\varepsilon, \varepsilon' \in (0,1)$.  
\begin{align*}
 \|\Delta_{2m,v,z}f\|^2 &\geq (1-\varepsilon)\left\|\frac{1}{2m-\gamma v(Q)+\xi}
 S f\right\|^2 + \left(1- \frac{1}{\varepsilon}\right)
 \left\|\frac{\gamma (\alpha^+\cdot  \nabla v)(Q)}{(2m-\gamma
   v(Q)+\xi)^2} \alpha^- \cdot P  f\right\|^2.  
\\
&\geq (1- \varepsilon )\frac{1}{1+ 16m^2} \left\|S f\right\|^2 
+ \left(1- \frac{1}{\varepsilon}\right) \frac{ \kappa \|\alpha^+\cdot
  \nabla v(Q) \| }{\delta^4} \left\|\alpha^- \cdot P f\right\|^2,
\\
&\hspace*{-2cm}\geq \left((1- \varepsilon )\frac{1}{1+ 16m^2} + \varepsilon' 
(\varepsilon-1) \frac{ \kappa \|\alpha^+\cdot
  \nabla v(Q) \| }{2 \varepsilon \delta^4}\right) \left\|S f\right\|^2 
+ (\varepsilon- 1) \frac{ \kappa \|\alpha^+\cdot
  \nabla v(Q) \| }{2 \varepsilon\varepsilon' \delta^4}
\left\|f\right\|^2. 
\end{align*} 
Choosing $ \varepsilon'$ small enough, we infer
\eqref{e:mourrestrict0chk}. 

We turn to \eqref{e:commu21}. Again, it is enough to compute in the
form sense on  $\Cc^\infty_c(\R^3;\C^\nu)$. 
\begin{align*}
[[\Delta_{2m,v,z},\rmi A],\rmi A]
=&\, 4\, \alpha^+\cdot P\frac{1}{2m -v +z}\alpha^-\cdot P
+4\, \alpha^+\cdot P\frac{Q\cdot \nabla v(Q)}{(2m -v +z)^2}\alpha^-\cdot P\\ 
&- \alpha^+\cdot P\frac{Q\cdot \nabla\big(Q\cdot \nabla v(Q)\big)}{(2m
  -v +z)^2}\alpha^-\cdot P+ 2\, \alpha^+\cdot P\frac{(Q\cdot \nabla
  v(Q))^2}{(2m -v +z)^3}\alpha^-\cdot P + (Q\cdot \nabla)^2 v(Q).
\end{align*}
Note that (H1) ensures that $\|(Q\cdot \nabla)^2 v(Q) f\|^2\leq 4 \|\,
|Q|(Q\cdot \nabla)^2 v(Q)\|^2 \|S f\|^2$ is controlled by $S$. 
Relying again on \eqref{e:bd}, the bound \eqref{e:commu21} follows. \qed

We finally turn to the proof of the main result of this section.
\proof[Proof of Theorem \ref{t:mourreD}] Since we have $e^{itA}\Hr^2\subset
\Hr^2$ and \eqref{e:commuS}, we obtain that
$\Delta_{2m,\gamma v,\xi}\in\Cc^1(A, \Hr^2, \Hr)$, for all $(\gamma,
\xi)\in \Ec$. By interpolation, we deduce that  $\Delta_{2m,\gamma
  v,\xi}\in\Cc^1(A, \Hr^1, \Hr^{-1})$. Now taking in account
\eqref{e:commu21}, we infer $\Delta_{2m,\gamma
  v,\xi}\in\Cc^2(A, \Hr^1, \Hr^{-1})$, for all $(\gamma,
\xi)\in \Ec$. 

Using Propositions \ref{p:commu1} and \ref{p:commu2}, we can apply
Theorem \ref{t:mourrestrict}. We derive there is a finite $C'$ so that
\begin{eqnarray*}
 \sup_{\Re z\geq 0,\Im z>0, (\gamma, \xi)\in \Ec}
|\langle f, (\Delta_{2m,\gamma v,\xi}-z)^{-1} f\rangle | \leq C'
\left(\|S^{-1/2} f\|^2 + \|S^{-1/2} A f\|^2\right).
\end{eqnarray*} 
The Hardy inequality concludes. \qed

\section{Main result}\label{s:mainresult}
In this section, we will prove the main result of this paper and
deduce Theorem \ref{t:main}. 

\begin{theorem}\label{t:vraimain}
Let $\gamma\in \R$. Suppose that $v\in L^\infty(\R^3; \R)$ satisfies
the hypothesis: 
\begin{enumerate}
\item[(H1)] $\|v\|_\infty\leq m/2$ and $\nabla v$,  $ Q\cdot \nabla v(Q)$,
 $\langle Q\rangle (Q\cdot \nabla v)^2(Q)$ are bounded. 
\item[(H2)] There are $c_v\in [0,2)$ and $c_v'\geq 0$ such that
\begin{eqnarray*}
 x\cdot
  (\nabla v)(x) + c_v v(x) \leq \frac{c_v'}{|x|^2}, \mbox{ for all } x\in \R^3\setminus\{0\}.
\end{eqnarray*} 
\end{enumerate}
Set $V_1(Q):= v(Q)\otimes \Id_{\C^{2n}}$, where $L^2(\R^3; \C^{2\nu})\simeq
L^2(\R^3)\otimes \C^{2\nu}$.
\begin{enumerate}
\item[(H3)] Consider $V_2\in L^1_{\rm loc}(\R^3; \R^{2\nu})$ satisfying: 
\begin{eqnarray*}
\langle Q\rangle^2 V_2(Q)  \in \Bc\big(\Hr^{1}(\R^3; \C^{2\nu}), L^2
(\R^3; \C^{2\nu})\big). 
\end{eqnarray*} 
\end{enumerate}
Then, there are $\kappa, \delta, C>0$, such that $H_\gamma:= D_m + \gamma V(Q)$,  where $V:=V_1 + V_2$, is self-adjoint with domain $\Hr^1(\R^3; \C^{2\nu})$. Moreover,
\begin{eqnarray}\label{e:vraimain}
\sup_{|\lambda|\in [m, m+\delta], \, \varepsilon>0, |\gamma|\leq \kappa }\|\langle
Q\rangle^{-1}(H_\gamma -\lambda - \rmi \varepsilon )^{-1} \langle
Q\rangle^{-1}\|\leq C.	 
\end{eqnarray}	
In particular, $H_\gamma$ has no eigenvalue in $\pm m$. Moreover,
there is $C'$ so that 
\begin{eqnarray}\label{e:vraiKato}
\sup_{|\gamma|\leq \kappa}\int_\R \| \langle Q \rangle^{-1}
e^{-it H_\gamma } E_\Ic(H_\gamma) f \|^2 dt  \leq C'\|f\|^2,
\end{eqnarray} 
where $\Ic=[-m-\delta, -m]\cup [m, m+\delta]$ and where
$E_\Ic(H_\gamma)$ denotes the spectral measure of $H_\gamma$.
\end{theorem}

\begin{remark} 
In \cite{FournaisSkibsted} and in \cite{Richard}, 
one takes advantage that the Virial of the potential is negative, 
in order to prove the limiting absorption principle for some
self-adjoint Schr\"odinger operators, see Remark \ref{r:virial}. 
Here, we cannot allow this hypothesis as we are also interested in
positronic threshold, i.e., we seek a result for $v$ and $-v$, see
\eqref{e:gauge}. We recover the positivity using some Hardy inequality
and small coupling constants.
\end{remark} 

\proof[Proof of Theorem \ref{t:vraimain}] First note that
\eqref{e:vraiKato} is a consequence of 
\eqref{e:vraimain}, see \cite{Kato}. Consider the case
$V_2=0$. Note that, in Section \ref{s:red}, the operator $H_\gamma$ is
denoted by $H^{\rm}_\gamma$. 

The self-adjointness is clear. We first apply
Theorem \ref{t:mourreD} and obtain \eqref{e:mourreD}. By
choosing $\xi=z$ and as $\| |Q| f\| \leq \| \langle Q\rangle f\|$, we
infer \eqref{e:LAPchange0}. In turn, it implies
\eqref{e:LAPchange1'}. Finally, using the unitary transformation
$\alpha_5$, \eqref{e:vraimain} follows from \eqref{e:gauge}. For a
general $V_2$, we use Proposition \ref{p:collage}. \qed 

It remains to explains how to add the singular part $V_2$ of the
potential by perturbing the limiting absorption principal. This is
somehow standard. Note that unlike \cite{JensenNenciu}, for
instance, we do not distinguish the nature of the
singularity at the threshold energy, as we work with small
coupling constants.

\begin{proposition}\label{p:collage}
Assume that Theorem \ref{t:vraimain} holds true for $V_2=0$. 
Take now $V_2$ satisfying (H3). Then there is $\kappa' \in (0,
\kappa]$, so that
\begin{eqnarray*}
 H_\gamma:= D_m + \gamma (V + V_2)(Q)
\end{eqnarray*} 
is self-adjoint with domain $\Hr^1(\R^3;
\C^{2\nu})$, for all $|\gamma|\leq \kappa$. Moreover, 
\begin{eqnarray*}
\sup_{\Re z\in [m, m+\delta], \Im z>0, |\gamma|\leq \kappa'}
\left\|\langle Q\rangle^{-1}(H_\gamma -z)^{-1}\langle
Q\rangle^{-1}\right\|<\infty. 
\end{eqnarray*} 
\end{proposition} 
\proof Up to a smaller $\kappa$, Kato-Rellich ensures the
self-adjointness. We turn to the estimate of the
resolvent. Easily, one reduces to the case $|\Im(z)|\leq 1$. From the
resolvent identity, we have: 
\begin{align*}
\nonumber
\langle Q\rangle^{-1} (H_\gamma -z)^{-1}\langle Q\rangle^{-1} \quad \langle
Q\rangle \left\{ I+\gamma V_2(H_\gamma^{\rm bd}
-z)^{-1}\right\}\langle  Q\rangle^{-1} 
=&\, \langle
Q\rangle^{-1}  (H_\gamma^{\rm bd} -z)^{-1} \langle
Q\rangle^{-1}
\end{align*}
Considering Lemma \ref{l:LAPchange} and Theorem \ref{t:mourreD}, 
the result follows if we can invert the second term of the
l.h.s.\ uniformly in the parameters. Therefore, we show there is $\kappa'\in
(0,\kappa]$ so that 
\begin{eqnarray*}
\sup_{\Re(z)\in [m, m+\delta], \Im(z)\in (0,1], |\gamma|\leq \kappa'} \|\langle Q\rangle
  \gamma V_2(H_\gamma^{\rm bd}  -z)^{-1} \langle  Q\rangle^{-1}\| <1. 
\end{eqnarray*} 
Using the identity of the resolvent, we get
\begin{align*}
 \langle Q\rangle V_2(H_\gamma^{\rm bd} -z)^{-1}\langle Q\rangle^{-1}
 =& \langle Q\rangle V_2 (H_0^{\rm bd} -i)^{-1} \quad  \langle Q
 \rangle^{-1}
\\
&- \langle Q\rangle V_2 (H_0^{\rm bd} -i
)^{-1}\langle Q\rangle \quad 
(\gamma V -z+i) \quad \langle Q\rangle^{-1}  (H_\gamma^{\rm bd} -z )^{-1}
\langle Q\rangle^{-1}. 
\end{align*}
The first term of the r.h.s.\ is bounded by using (H3).  To
  control the last term, remember that $z$ is 
bounded and use again Lemma \ref{l:LAPchange} and Theorem
\ref{t:mourreD}. It remains to notice that 
\[\langle Q\rangle V_2  (H_0^{\rm bd} -i
)^{-1}\langle Q\rangle=\langle Q\rangle^2 V_2  (H_0^{\rm bd} -i
)^{-1}-\langle Q\rangle V_2  (H_0^{\rm bd} -i
)^{-1} \quad [H_0^{\rm bd},\langle Q\rangle]_\circ(H_0^{\rm bd} -i
)^{-1}\]
is bounded. Indeed, the assumption (H3) controls the terms in
$V_2$ and $\langle Q\rangle\in \Cc^1(H_0)$ and $[H_0^{\rm bd},\langle
  Q\rangle]_\circ$ is bounded, see proof of Lemma \ref{l:steptoLAP}.  \qed  

At last, Theorem \ref{t:main} is  an immediate corollary of Theorem 
\ref{t:vraimain}. Indeed, one has:

\begin{example}[Multi-center]\label{ex:multi} For $i=1, \ldots,
  n$, we choose 
  $a_i\in \R^3$ the site of the poles and $Z_i\in \R$ its charge. We
  set: 
\begin{eqnarray*}
v_c := \sum_{i=1}^n \frac{z_i}{|\cdot- a_i|}
\end{eqnarray*}
Note that 
\begin{eqnarray*}
Q\cdot \nabla v_c (Q) + v_c:= \sum_{i=1}^n a_i \cdot \nabla v(Q).
\end{eqnarray*}
Choose now $\varphi\in\Cc^\infty_c(\R^3)$ radial with values in
$[0,1]$. Moreover, we ask that $\varphi$ restricted to the ball $B(0,
\max(|a_i|))$ is $1$. Consider the support large enough, so that
$\|\tilde \varphi v\|_\infty \leq m/2$, where $\tilde
\varphi:\ 1-\varphi$. Set $v:= \tilde 
\varphi v_c$. Straightforwardly, the hypothesis (H1) and (H2) are
satisfied. Note that (H3) follows from the Hardy inequality. 
\end{example}

\begin{example}[Smooth homogeneous potentials] In \cite{Herbst}, one considers smooth potential independent of $|x|$ of the form $v(x):= \tilde v(|x|/x)$, with 
$v\in\Cc^\infty(S^{2})$, see also Remark \ref{r:Herbst}. Here, by taking $c_v=0$ in Theorems \ref{t:vraimain} and \ref{t:vraimain2}, one obtains a relativistic equivalent of this result. We point out that this perturbation is not relatively compact with respect to the Dirac operator. 
\end{example}

We now discuss singular weights in $|Q|$.

\begin{remark}\label{r:noLAP} It is important to note that unlike in the non-relativistic case, see Theorem \ref{t:nonrelamain}, one cannot replace the weights $\langle Q\rangle$ in \eqref{e:mainlow} by $|Q|$. 
Indeed, with the notation of Theorem \ref{t:vraimain}, $V_2=0$ and $z\in \C$, 
consider a function $f$ in $\Cc^\infty_c(\R^3\setminus \{0\};\C^{2\nu})$ and notice the expression $R^{-3/2}\big\|\,|Q|(H_\gamma-z) |Q| f(\cdot/R)\big\|_2$ tends to $0$, as $R$ goes to $0$. Therefore, there is no $z\in\C$ such that the operator $|Q|(H_\gamma -z) |Q|$ has a bounded inverse.
\end{remark}  

We finally give a second result with a weight allowing some singularity in $|Q|$.
Using Lemma \ref{l:LAPchange2} instead of Lemma \ref{l:LAPchange} in the proof of Theorem \ref{t:vraimain}, we infer straightforwardly:
\begin{theorem}\label{t:vraimain2}
Let $\gamma\in \R$ and take $v\in L^\infty(\R^3; \R)$ satisfying (H1) and (H2). 
Then, there are $\kappa, \delta, C>0$, such that 
$H_\gamma:= D_m + \gamma v(Q)\otimes \Id_{\C^{2\nu}}$ satisfies
\begin{eqnarray}\label{e:vraimain2}
\sup_{|\lambda|\in [m, m+\delta], \, \varepsilon>0, |\gamma|\leq \kappa }\|  \langle P\rangle^{-1/2} |Q|^{-1}(H_\gamma -\lambda - \rmi \varepsilon )^{-1} |Q|^{-1}\langle P\rangle^{-1/2}\|\leq C.	 
\end{eqnarray}	
Moreover,
there is $C'$ so that 
\begin{eqnarray}\label{e:vraiKato2}
\sup_{|\gamma|\leq \kappa}\int_\R \| \langle P\rangle^{-1/2} |Q|^{-1}
e^{-it H_\gamma } E_\Ic(H_\gamma) f \|^2 dt  \leq C'\|f\|^2,
\end{eqnarray} 
where $\Ic=[-m-\delta, -m]\cup [m, m+\delta]$ and where
$E_\Ic(H_\gamma)$ denotes the spectral measure of $H_\gamma$.
\end{theorem}

Keeping in mind Proposition \ref{p:collage}, one sees that one can
only add trivial potentials in the perturbation theory of the limiting
absorption principle. Thence, it is an open question whether one can
cover the example \ref{ex:multi} with the weights $\langle P\rangle^{1/2}|Q|$.

\appendix
\section{Commutator expansions.} \label{s:dev-commut}
\setcounter{equation}{0}  
This section is a small improvement of \cite{GoleniaJecko}[Appendix B], see also \cite{DerezinskiGerard, HunzikerSigalSoffer}. 
We start with some generalities. Given a bounded operator $B$ and a
self-adjoint operator $A$ acting in a Hilbert space $\Hr$, one says
that $B\in \Cc^k(A)$ if $t\mapsto e^{-itA}B e^{itA}$ is strongly
$\Cc^k$. Given a self-adjoint operator $B$, one says that $B\in
\Cc^k(A)$ if for some (hence any) $z\notin \sigma(B)$, $t\mapsto
e^{-itA}(B-z)^{-1} e^{itA}$ is strongly $\Cc^k$. The two definitions
coincide in the case of a bounded self-adjoint operator. We recall a
result following from Lemma 6.2.9 
and Theorem 6.2.10 of \cite{AmreinBoutetdeMonvelGeorgescu}. 
\begin{theorem}\label{th:abg} 
Let $A$ and $B$ be two self-adjoint operators in the Hilbert space
$\Hr$. For $z\notin \sigma(A)$, set $R(z):=(B-z)^{-1}$. The following
points are equivalent to $B\in\Cc^1(A)$:  
\begin{enumerate} 
\item For one (then for all) $z\notin \sigma(B)$, there is a finite
$c$ such that 
\begin{align}\label{e:C1a} 
|\langle A f, R(z) f\rangle - \langle R(\overline{z}) f, Af\rangle| \leq c 
\|f\|^2, \mbox{ for all $f\in\Dc(A)$}. 
\end{align} 
\item 
\begin{enumerate} 
\item [a.]There is a finite $c$ such that for all $f\in \Dc(A)\cap\Dc(B)$: 
\begin{equation}\label{e:C1b} 
|\langle Af, B f\rangle- \langle B f, Af\rangle|\leq \,
 c\big(\|B f\|^2+\|f\|^2\big). 
\end{equation} 
\item [b.] For some (then for all) $z\notin \sigma(B)$, the set
$\{f\in\Dc(A) \mid R(z)f\in\Dc(A)$ and $R(\overline{z})f\in\Dc(A)
\}$ is a core for $A$. 
\end{enumerate} 
\end{enumerate} 
\end{theorem} 
Note that the condition (3.b) could be uneasy to check, see
\cite{GeorgescuGerard}. We mention \cite{GoleniaMoroianu}[Lemma A.2]
and \cite{GerardLaba}[Lemma 3.2.2] to overcome this subtlety. 
As $(B+i)^{-1}$ is a homeomorphism between $\Hr$ onto $\Dc(B)$, 
$(B+i)^{-1} \Dc(A)$ is dense in $\Dc(B)$, endowed with the graph norm. 
Moreover, \eqref{e:C1a} gives $(B+i)^{-1}\Dc(A)\subset
\Dc(A)$. Therefore $(B+i)^{-1}\Dc(A)\subset \Dc(B)\cap \Dc(A)$ are
dense in $\Dc(B)$ for the graph norm. Remark that $\Dc(B)\cap \Dc(A)$ is usually not dense in $\Dc(A)$, see \cite{GeorgescuGerardMoller}. 

Note that \eqref{e:C1a} yields the commutator $[A, R(z)]$ extends to a
bounded operator, in the form sense. We shall denote the extension by
$[A, R(z)]_\circ$. In the same way, since $\Dc(B)\cap \Dc(A)$ is dense in $\Dc(B)$, 
\eqref{e:C1b} ensures that the commutator $[B, A]$ extends to a unique element of $\Bc\big(\Dc(B), \Dc(B)^*\big)$ denoted by $[B,   A]_\circ$. Moreover, when $B\in
\Cc^1(A)$, one has: 
\begin{eqnarray*}
\big[A, (B-z)^{-1}\big]_\circ =\quad  \underbrace{(B-z)^{-1}}_{\Hr
  \leftarrow \Dc(B)^*}\quad  \underbrace{[B, A]_\circ}_{\Dc(B)^*\leftarrow
  \Dc(B)} \quad \underbrace{(B-z)^{-1}}_{\Dc(B)\leftarrow \Hr}.
\end{eqnarray*} 
Here we use the Riesz lemma to identify $\Hr$ with its anti-dual
$\Hr^*$. 

We now recall some well known facts on symbolic calculus and almost
analytic extensions. For $\rho\in\R$, let $\Sc^\rho$ be the class of
function $\varphi\in\Cc^\infty(\R;\C)$ such that    
\begin{eqnarray}\label{eq:regu} 
\forall k\in\N, \quad C_k(\varphi) :=\sup _{t\in\R}\, \langle
t\rangle^{-\rho+k}|\varphi^{(k)}(t)|<\infty .      
\end{eqnarray} 
Equipped with the semi-norms defined by (\ref{eq:regu}), $\Sc^\rho$ is
a Fr\'echet space. Leibniz' formula implies the continuous embedding:
$\Sc^\rho\cdot  \Sc^{\rho'} \subset  \Sc^{\rho+\rho'}$.  
We shall use the following result, e.g., \cite{DerezinskiGerard}. 
 
\begin{lemma}\label{l:dg} 
Let $\varphi\in\Sc^\rho$ with $\rho\in\R$. For all   
$l\in \N$, there is a smooth function  $\varphi^\C:\C \rightarrow \C$,
such that: 
\begin{eqnarray} 
\label{eq:dg1} \varphi^\C|_{\R}=\varphi,\quad &&\left|\frac{\partial 
  \varphi^\C}{\partial  \overline{z}}(z) \right|\leq c_1 \langle \Re(z) 
  \rangle^{\rho-1 -l} |\Im(z)|^l\\\label{eq:dg2} 
&& \supp \varphi^\C \subset\{x+iy\mid |y|\leq c_2 \langle 
  x\rangle\},\\\label{eq:dg3} && \varphi^\C(x+iy)= 0, \mbox{ if } 
  x\not\in\supp\varphi .  
\end{eqnarray}  
for some constants $c_1$, $c_2$ depending on the semi-norms 
\eqref{eq:regu} of $\varphi$ in $\Sc^\rho$ and not on $\varphi$.   
\end{lemma}  
One calls $\varphi^\C$ an \emph{almost analytic extension} of $\varphi$.
Let $A$ be a self-adjoint operator, $\rho < 0$ and $\varphi\in
\Sc^{\rho}$. By functional calculus, one has $\varphi(A)$
bounded. The Helffer-Sj\"{o}strand's formula, see
\cite{HelfferSjostrand} and \cite{DerezinskiGerard} for instance,
gives that for all almost analytic extension of $\varphi$, one has:  
\begin{eqnarray}\label{eq:int} 
\varphi(A) = \frac{i}{2\pi}\int_\C\frac{\partial \varphi^\C }{\partial 
\overline{z}}(z-A)^{-1}dz\wedge d\overline{z}. 
\end{eqnarray}  
Note the integral exists in the norm topology, by \eqref{eq:dg1} 
with $l=1$. Next we come to a commutator expansion. Here $B$ is not
necessarily bounded while in \cite{GoleniaJecko}, one 
considers the case $B$ bounded. We denote by $\ad_A^j(B)$ the
extension of the $j$-th commutator of $A$ with $B$ defined inductively
by $\ad_A^p(B):=[\ad_A^{p-1}(B),A]_\circ$, when it exists. 
 
\begin{proposition}\label{p:regu} 
Let $k\in \N^\ast$ and $B\in\Cc^k(A)$ be self-adjoint. Suppose 
$\ad_A^j(B)$ are bounded operators, for $j=1,\ldots, k$. Let 
$\rho < k$ and $\varphi\in \Sc^{\rho}$. Suppose that $\Dc(B)\cap \Dc(\varphi(A))$ is dense in $\Dc(\varphi(A))$ for the graph norm. Then, the commutator
$[\varphi(A), B]_\circ$ belongs to  $\Bc\big(\Dc(\varphi'(A)), \Hr\big)$ and satisfies 
\begin{eqnarray}\label{e:ega}
[\varphi(A), B]_\circ = \sum_{j=1}^{k-1} \frac{1}{j!} 
\varphi^{(j)}(A)\ad_A^j(B)   
+ \frac{i}{2\pi}\int_\C\frac{\partial\varphi^\C }{\partial 
\overline{z}}(z-A)^{-k} \ad_A^k(B) (z-A)^{-1} dz\wedge d\overline{z}, 
\end{eqnarray} 
where the integral exists for the topology of $\Bc(\Hr)$. 
\end{proposition} 
\proof We cannot use the \eqref{eq:int} directly with $\varphi$ as the
integral does not seem to exist. We proceed as in
\cite{GoleniaJecko}. Take $\cchi_1\in \Cc^\infty_c(\R;\R)$ with values
in $[0,1]$ and being $1$ on $[-1,1]$. Set
$\cchi_R:=\cchi(\cdot/R)$. As $R$ goes to infinity, $\cchi_R$
converges pointwise to $1$. Moreover, $\{\cchi_R\}_{R\in [1,\infty]}$ is bounded in $\Sc^0$. We infer $\varphi_R:=\varphi \cchi_R$ tends pointwise 
to $\varphi$ and that $\{\varphi_R\}_{R\in [1,\infty]}$ is bounded in $\Sc^\rho$. Now, note that
\begin{align}\label{e:ligne}  
\hspace*{1cm} 
[\varphi_R(A), B] &  
=\sum_{j=1}^{k-1} \frac{i}{2\pi}\int_\C\frac{\partial\varphi^\C_R 
}{\partial\overline{z}}  (z-A)^{-j-1}\ad_A^{j}(B) dz\wedge 
d\overline{z}
\\ 
\nonumber
&\, +  \frac{i}{2\pi}\int_\C\frac{\partial\varphi^\C_R }{\partial 
\overline{z}}(z-A)^{-k} \ad_A^k(B) (z-A)^{-1} dz\wedge d\overline{z}. 
\end{align} 
in the form sense on $\Dc(B)$. Using \eqref{eq:dg1}, the integral
converges in norm. We write $[\varphi_R(A), B]_\circ$ on the
l.h.s. The first term of the r.h.s.\ is  $\sum_{j=1}^{k-1} 
\varphi^{(j)}_R(A)\ad_A^j(B)/ j!$. Now we let $R$ goes to infinity. On the l.h.s.\ we use the Lebesgue converges. On the r.h.s.\ we expand the commutator in \eqref{e:ligne} in the form sense on $\Dc\big(\varphi(A)\big)\cap \Dc(B)$, take the limit by functional calculus and finish by density in $\Dc\big(\varphi(A)\big)$. \qed 

The hypothesis on the density of $\Dc(B)\cap \Dc(\varphi(A))$ in $\Dc(\varphi(A))$ could be delicate to check. It follows by the Nelson Lemma from the fact that the $C_0$-group $\{e^{it A^k}\}_{t\in \R}$ stabilizes $\Dc(B)$. We mention that for $k=1$, since $[B,iA]_\circ$ is bounded, \cite{GeorgescuGerard}[Lemma 2] ensures this invariance of the domain.

The rest of the previous expansion is estimated as in
\cite{GoleniaJecko}. We rely on the following important bound. Let
$c>0$ and $s\in [0,1]$, there exists some $C>0$ so that, for all
$z=x+iy\in\{a+ib\mid 
0<|b|\leq c\langle a\rangle \}$:  
\begin{eqnarray}\label{eq:majoA} 
\big\| \langle A\rangle^s (A-z)^{-1}\big\|\leq C \langle x
\rangle^{s}\cdot |y|^{-1}.  
\end{eqnarray}

\begin{lemma}\label{l:est3} 
Let $B\in\Cc^k(A)$ self-adjoint. Suppose 
$\ad_A^j(B)$ are bounded operators, for $j=1,\ldots, k$. 
 Let $\varphi\in\Sc^\rho$, with $\rho< k$. Let $I_k(\varphi)$ the
 rest of the development of order $k$ of $[\varphi(A), B]$ in
 \eqref{e:ega}. Let  $s, s' \in [0,1]$ such that 
 $\rho+s+s'<k$. Then $\langle A \rangle^{s}  I_k(\varphi)\langle A
 \rangle^{s'}$ is bounded and it is uniformly bounded when $\varphi$
 stays in a bounded subset of $\Sc^\rho$. Let $R>0$. If $\varphi$ stays
 in a bounded subset of $\{\psi \in \Sc^\rho\mid [-R;R]\cap
 \supp(\varphi)=\emptyset\}$ then $\langle R\rangle^{k-\rho-s-s'}
 \|\langle A \rangle^{s} I_k(\varphi)\langle A \rangle^{s'}\|$ is
 uniformly bounded.   
\end{lemma}  
\proof 
We will follow ideas from \cite{DerezinskiGerard}[Lemma C.3.1]. In
this proof, all the  constants are denoted by $C$, independently of
their value. Given a complex  number $z$, $x$ and $y$ will denote its
real and imaginary part, respectively. Since $B\in\Cc^k(A)$,
$\ad^k_A(B)$ is bounded.  We start with the second assertion. Let
$\varphi\in \Sc^\rho$, $R>0$ such that $[-R;R]\cap
\supp(\varphi)=\emptyset$. Notice that,  by \eqref{eq:dg3},
$\varphi^\C(x+iy)= 0$ for $|x|\leq R$. By \eqref{eq:majoA}, 
\begin{align*} 
\|\langle A \rangle^{s}I_k(\varphi)\langle A \rangle^{s'}\|\leq&\,
\frac{1}{\pi} \int  
\big|\frac{\partial\varphi^\C }{\partial \overline{z}}\big|\cdot 
\frac{\langle   x\rangle^{s}}{|y|^k}  \cdot  
\|\ad^k_A(B)\|\cdot \frac{\langle x\rangle^{s'}}{|y|} dx\wedge dy\\ 
\leq&\,  C(\varphi)\int_{|x|\geq R}\int_{|y|\leq c_2\langle x\rangle } \langle x 
\rangle^{\rho+s+s'-1-l}|y|^l |y|^{-k-1} dx\wedge dy,  
\end{align*}  
for any $l$, by \eqref{eq:dg1}. Recall that $dz\wedge d\overline{z}=-2i dx\wedge  dy$. We choose $l=k+1$. We have,  
\begin{eqnarray*} 
\|\langle A\rangle^{s}I_k(\varphi)\langle A\rangle^{s'} 
\|\leq C(\varphi)\int_{|x|\geq R} \langle x \rangle^{\rho+s+s'-k-1}dx 
\leq C(\varphi)\langle R\rangle^{\rho+s+s'-k}.   
\end{eqnarray*}  
Since $C(\varphi)$ is bounded when $\varphi$ stays in a bounded subset
of $\Sc^\rho$, this yields the second assertion. For the first one, we can follow the same lines, 
replacing $R$ by $0$ in the integrals, and arrive at the result. \qed 

\section{A non-selfadjoint weak Mourre theory}\label{s:WeakMourreTheory}
In this section, we adapt ideas coming from \cite{FournaisSkibsted}
and \cite{Richard} in order to obtain a limiting absorption
principle for a family of closed operators $\{H^{\pm}(p)\}_{p\in \Ec}$. 
We ask that they have a common domain 
\begin{eqnarray}\label{e:dom0}
\Dr:=\Dc\big(H^+(p)\big)=\Dc\big(H^-(p)\big), \mbox{  for all } p\in \Ec.
\end{eqnarray}
We choose $p_0\in \Ec$ and endow $\Dr$ with the graph norm of
$H^+(p_0)$. We also ask that
\begin{eqnarray}\label{e:adj}
\big(H^+(p)\big)^*=H^-(p), \mbox{  for all } p\in \Ec.
\end{eqnarray}
In particular, we have that $\Dc\big((H^{\pm}(p))^*\big)= \Dr$. In the
sequel, we forgo $p$, when no confusion can arises.

Since $H^\pm$ are densely defined, share the same domain and are
adjoint of the other, we have that $\Re (H^\pm)$ and $\Im (H^\pm)$ are
closable operators on $\Dr$, indeed their adjoints are densely
defined. We denote by  $\Re (H^\pm)$ and by   $\Im (H^\pm)$ the
closure of these operators. It is possible that they are not
self-adjoint, albeit there are symmetric. However, $\Dr$ is a core for
them. Their domain is possibly bigger than $\Dr$. 
We suppose that $H^+$ is \emph{dissipative}, i.e.,
\begin{eqnarray*}
\mbox{  $\langle f, \Im(H^+) f \rangle  \geq 0$, for all $f\in \Dr$.} 
\end{eqnarray*} 
This gives also that $\Im (H^-)\leq 0$. By the numerical range
theorem (see Lemma \ref{l:NRT}), we infer that $\sigma(H^{\pm})$ is included in the half-plan containing $\pm i$.
Take now a non-negative self-adjoint operator $S$,
\emph{independent} of $p\in \Ec$, with form domain
$\Gr:=\Dc(S^{1/2})\supset\Dr$.  We assume that $S$ is injective. We have $\langle f, S  f\rangle> 0$  for
all $f\in \Gr\setminus \{0\}$ and simply write $S> 0$. One defines
$\Sr$ as the completion of $\Gr$ under the norm  $\|f\|_\Sr^2:=\langle
f, S f\rangle$. We obtain $\Gr\subset \Sr$ with dense and continuous
embedding. Moreover, since $\Gr= \langle S^{1/2}\rangle^{-1}\Hr$,
$\Sr$ is also the completion of $\Hr$ under the norm   given by
$\|S^{1/2} \langle S^{1/2}\rangle^{-1}\cdot\|$. We use the Riesz Lemma
to identify $\Hr$ with $\Hr^*$, its anti-dual. The adjoint space $\Sr^*$ of $\Sr$ is
exactly the domain of $\langle S^{1/2}\rangle S^{-1/2}$ in
$\Hr\simeq\Hr^*$. Note that $S^{-1}$ is an isomorphism between  $\Sr$ and
$\Sr^*$.  We get the following scale with continuous and dense embeddings:    
\begin{eqnarray}\label{e:scale}
\begin{array}{cccccccccc}
&&&&&& \Sr^*&&&
\\
&&&&&& \downarrow&\searrow &&
\\
\Dr&\longrightarrow &\Gr& \longrightarrow &\Hr&\simeq& \Hr^*
&\longrightarrow &\Gr^*& \longrightarrow\Dr^*.
\\
&&&\searrow &\downarrow &&&&&
\\
&&&&\Sr&&&&&
\end{array}
\end{eqnarray}   
To perform this analysis, we consider an external operator, the
conjugate operator. Let $A$ be a self-adjoint operator in $\Hr$. 
We assume $S\in \Cc^1(A)$. Let $W_t:=e^{itA}$ be the $C_0$-group
associated to $A$ in $\Hr$. We ask:
\begin{eqnarray}\label{e:stab}
 W_t\Gr\subset \Gr \mbox{ and } W_t \Sr\subset \Sr, \mbox{ for all
 } t\in\R. 
\end{eqnarray} 
By duality, we have $W_t$ stabilizes $\Gr^*$ and also $\Sr^*$ (but may
be not $\Dr$ or $\Dr^*$). The restricted group to these spaces is also
a $C_0$-group.  We denote the generator by $A$ with the subspace in
subscript.  Given $\Hr_i\subset\Hr_j$ be two of those spaces. One
easily shows that $A|_{\Hr_i}\subset A|_{\Hr_j}$ and that
$A|_{\Hr_j}$ is the closure of $A|_{\Hr_i}$ in $\Hr_j$. Moreover, one has
\begin{eqnarray}\label{e:dom}
\Dc(A|_{\Hr_i})= \left\{f\in \Dc\big(A|_{\Hr_j}\big)\cap\Hr_i \mbox{
such that } A|_{\Hr_j} f\in \Hr_i\right\}. 
\end{eqnarray}
We now explain how to check the second hypothesis of \eqref{e:stab},
see also \cite{Richard}.  
\begin{remark}\label{r:inv}
The second invariance of the domains of \eqref{e:stab} follows from
the first one  and from
\begin{eqnarray}
|\langle S f, A f\rangle- \langle A f, S f\rangle| \leq c \|S^{1/2}
f\|^2, \mbox{ for all } f\in\Dc(S)\cap\Dc(A). 
\end{eqnarray} 
As $(S+i)^{-1}$ is a homeomorphism between $\Hr$ onto $\Dc(S)$, 
$(S+i)^{-1} \Dc(A)$ is dense in $\Dc(S)$, endowed with the graph norm. 
Moreover, since $S\in \Cc^1(A)$, one has $(S+i)^{-1}\Dc(A)\subset
\Dc(A)$. Therefore $(S+i)^{-1}\Dc(A)\subset \Dc(S)\cap \Dc(A)$ are
dense in $\Dc(S)$, hence in $\Gr$ and in $\Sr$. 
The commutator $[S,A]$ has a unique extension to an element of $\Bc(\Sr, \Sr^*)$, in
the form sense. We denote it by $[S,A]_\circ$. Take now $f\in
\Gr\cap\Dc(A)$, which is a dense set in $\Gr$. On one hand we have
$\tau\mapsto \|W_\tau f\|_\Sr^2$ is bounded when $\tau$ is in a
compact set (since $\Gr \, \hookrightarrow \Sr$. On
the other  hand, the Gronwall lemma concludes by noticing:   
\begin{eqnarray*}
\|W_t f\|_\Sr^2	= \langle f, S f\rangle + \int_0^t \langle W_\tau f,
   [S,iA]_\circ W_\tau f \rangle\, d\tau 
\leq \|S^{1/2} f\|^2 + c \int_0^{|t|} \|W_\tau f\|_\Sr^2\, d\tau.  
\end{eqnarray*}
\end{remark}

Let $\Kr\subset \Hr$ be a space which is stabilized by $W_t$. 
Consider $L\in \Bc(\Kr, \Kr^*)$. We say that $L\in \Cc^k(A; \Kr,
\Kr^*)$, when $t\rightarrow W_{-t}L W_t$ is strongly $\Cc^k$ 
from $\Kr$ into $\Kr^*$. When $\Kr=\Hr$, this class is the same as 
$\Cc^k(A)$, see \cite{AmreinBoutetdeMonvelGeorgescu}[Theorem 6.3.4
a.]. 

\begin{theorem}\label{t:mourrestrict}
Let $H^{\pm}=H^{\pm}(p)$, with $p\in \Ec$ as above. Let  $A$ be self-adjoint 
such that \eqref{e:stab} holds true. Suppose that $H^\pm\in\Cc^2(A;
\Gr, \Gr^*)$ and that there is a constant $c$, independent of $p$, such that
\begin{eqnarray}\label{e:mourrestrict0}
|\langle H^{\mp}f, Ag \rangle - \langle Af, H^{\pm} g \rangle   |\leq 
c\|f\|\cdot\|(H^{\pm} \pm i) g\|, \mbox{ for all } f,g \in \Dr\cap \Dc(A).
\end{eqnarray}
Take $c_1\geq 0$ independent of $p$ and assume that   
\begin{eqnarray}
\label{e:mourrestrict1}	 [\Re (H^\pm), \rmi A]_\circ -c_1\Re (H^\pm) \geq S
	> 0,  &\\
\label{e:mourrestrict2}	\pm c_1[\Im (H^\pm),\rmi A]_\circ\geq 0, &
	\pm \Im (H^\pm)\geq 0,
\end{eqnarray}
in the sense of forms on $\Gr$. Suppose also there exists $C>0$
independent of $p\in \Ec$ such that 
\begin{equation}\label{Eq:SecondCommBound}
\left |\langle f, \big[\big[ H^\pm, A\big]_\circ,A\big]_\circ
f\rangle\right|\leq C  \|S^{1/2} f\|^2, \mbox{ for all } f\in\Gr.
\end{equation}
Then, there are $c$ and $\mu_0>0$, both independent of $p$,  such that
\begin{eqnarray}\label{e:LAP}
|\langle f, (H^\pm -\lambda \pm i \mu)^{-1} f\rangle | \leq c
\left(\|S^{-1/2} f\|^2 + \|S^{-1/2} A f\|^2\right)\leq c \|f\|_{\Dc(A|_{\Sr^*})},  
\end{eqnarray} 
for all $p\in \Ec$, $\mu\in (0, \mu_0)$ and $\lambda\geq 0$, in the
case $c_1>0$ and $\lambda\in \R$ if $c_1=0$.
\end{theorem}
In the self-adjoint setting, the case $c_1=0$ is treated in
\cite{BoutetKazantsevaMantoiu, BoutetMantoiu}. Comparing with
\cite{Richard}, who deal with the case of one 
self-adjoint operator and for $c_1>0$.  
We give some few improvements. First, we do not ask $\Dr$
to be the domain of $S$. Moreover, we drop the hypothesis
that the first commutator $[H, \rmi A]_\circ$ is bounded from
below. For the latter, we use more carefully the numerical range
theorem in our proof. Finally, unlike \cite{Richard}, we shall not go into
interpolation theory so as to improve the norm in the limiting
absorption principle. Indeed, in the context of the model we are
considering here, we reach the weights we are interested in without
it. We stick to an intermediate and explicit result, which is closer
to \cite{IftimoviciMantoiu}. Therefore, for the sake of clarity, we
present then the easiest proof possible and pay an important care
about domains.   

We also mention that there exists other Mourre-like theory for non-self-adjoint
operators, \cite{ABCF, Royer}. 

\proof We focus on the case $c_1>0$, as for the case $c_1=0$, one
replaces ``$\lambda>0$'' by ``$\lambda\in \R$''. We define $H_\varepsilon^\pm
:= H^\pm \pm\rmi \e[ H^\pm,\rmi A]_\circ$ with the common domain $\Dr$
for $\varepsilon\geq 0$.  Since $H^\pm \pm \rmi$ is bijective, by
writing $H_\varepsilon^\pm \pm \rmi= \big(1 \pm \rmi \e [ H^\pm,\rmi
  A]_\circ(H^\pm \pm \rmi)^{-1}\big)(H^\pm \pm \rmi)$ 
and using \eqref{e:mourrestrict0}, we get
there is $\varepsilon_0$ such that $H_\varepsilon^\pm(p) \pm \rmi$ is
bijective and closed for all $|\varepsilon|\leq
\varepsilon_0$ and all $p\in \Ec$. Therefore $(H_\varepsilon^\pm \pm \rmi)^*$ is also
bijective from $\Dc\big((H_{\varepsilon}^\pm)^*\big)$ onto $\Hr$. Now since
$(H_\varepsilon^\pm\pm\rmi)^*$ is an extension of
$H_{\varepsilon}^\mp \mp\rmi$ which is also bijective, we
infer the equality of the domains and that $(H_\varepsilon^\pm)^*=
H_{\varepsilon}^\mp$ for $\varepsilon\leq
\varepsilon_0$. 

Since $H^\pm\in \Cc^1(A; \Gr,\Gr^*)$, we obtain that $\Re(H^\pm)$ and
$\Im(H^\pm)$ are in $\Cc^1(A; \Gr,\Gr^*)$. In this space we have $[ H^\pm,
A]_\circ= [ \Re(H^\pm), A]_\circ + i [ \Im(H^\pm), A]_\circ$. Now,
take $f\in \Gr$.  Take $\varepsilon, \lambda, \mu \geq 0$.  We get: 
\begin{align}
\nonumber
&\hspace*{-1cm} -c_1 \varepsilon  \left\langle
f,\Re(H^\pm_\varepsilon-\lambda  \pm \rmi \mu   )f\right\rangle \mp
\left\langle f,\Im (H^\pm_\varepsilon-\lambda  \pm 
\rmi \mu ) f\right\rangle= 
\\ 
\nonumber
&=-c_1 \varepsilon \left\langle f,\left(\Re(H^\pm) \pm \e[ \Im
  (H^\pm),\rmi A]_\circ  -\lambda\right)f\right\rangle \mp
\left\langle f, \left(\Im( H^\pm)\mp\mu\mp\e[\Re (H^\pm),
  \rmi A]_\circ\right)f\right\rangle 
\\
\nonumber
&=\varepsilon \left\langle f,\big([\Re (H^\pm), \rmi A]_\circ -c_1\Re (H^\pm)
\big)f\right\rangle
+\left(c_1\lambda \varepsilon  +\mu \right)\left\|f\right\|^2
\mp\left\langle f, \left(c_1\e^2 [ \Im (H^\pm),\rmi A]_\circ 
+\Im (H^\pm)\right)f \right \rangle
\\
\label{e:numericalrange}
&\geq (c_1\lambda \varepsilon  + \mu)\left\|f\right\|^2 + \varepsilon
\| S^{1/2} f\|^2. 
\end{align}
We start with a crude bound. For $\varepsilon, \mu>0$, we get:
\begin{eqnarray*}
(c_1 \varepsilon +1)\, \|(H_\varepsilon^\pm-\lambda \pm\rmi\mu)f \|_{\Gr^*}\geq
 {\min(c_1\lambda \varepsilon  +\mu, \varepsilon ) } \|f \|_\Gr. 
\end{eqnarray*}
Since $H_\varepsilon^\pm-\lambda \pm\rmi\mu\in\Bc(\Gr, \Gr^*)$ and since they are adjoint
of the other, we infer the injectivity and that the ranges are
closed. They are bijective and the inverse is bounded by the open mapping theorem. 
\begin{eqnarray*}
 G^\pm_\e:=G^\pm_\e(\lambda,
\mu)=(H_\varepsilon^\pm-\lambda \pm\rmi\mu)^{-1} \mbox{ exists in }
\Bc(\Gr^*, \Gr), \mbox{ for } \lambda \geq 0 \mbox{
  and } \varepsilon, \mu>0. 
\end{eqnarray*} 
Here we lighten the notation but keep in mind the dependency in
$\lambda$  and $\mu$. Moreover, 
\begin{eqnarray}\label{e:G0}
\|G_\varepsilon^\pm \|_{\Bc(\Gr^*, \Gr)}\leq (c_1 \varepsilon
  +1)/\min(c_1\lambda \varepsilon  +\mu, \varepsilon), \mbox{ for }
  \lambda \geq 0 \mbox{ 
  and } \varepsilon, \mu>0. 
\end{eqnarray} 
This bound seems not enough to lead the whole analysis. Then, we first
restrict the domain of $G_\varepsilon^\pm$ to $\Hr$ and improve
it. Since this inequality 
\eqref{e:numericalrange} holds also true on the common domain of
$H_\varepsilon^\pm$ (and of its adjoint), we can apply the numerical
range theorem, Lemma \ref{l:NRT}. Since $S\geq 0$, we get the spectrum of
$H_\varepsilon^+ -\lambda +\rmi\mu$ is contained in the lower
half-plane delimited by the equation $y\leq -c_1 \varepsilon
x-\mu$.  
Hence, for $\varepsilon\in (0, \varepsilon_0]$ and $\mu>0$,
$H_\varepsilon^\pm -\lambda \pm\rmi\mu$ is bijective and by taking
$\varepsilon_0$ smaller, one has the distance from $0$ to the boundary
of the cone bigger than $\mu/2$. Then,
\begin{eqnarray}\label{e:G1}
\|G^\pm_\e \|_{\Bc(\Hr)}\leq 2/ \mu, \mbox{ for } \mu>0 \mbox{ and }
\varepsilon\in[0, \varepsilon_0]. 
\end{eqnarray} 
Note also that $(G^\pm_\e)^*= G^\mp_\e$. Take $\varepsilon, \mu>0$. We fix $f\in \Hr$ and set:
\begin{eqnarray*}
F_\e^\pm:=\left\langle f, G^\pm_\e f\right\rangle. 
\end{eqnarray*} 
Since $G^\pm_\e \Hr\subset \Dr\subset \Sr$ and using
\eqref{e:numericalrange},  we infer
\begin{align}
\nonumber
\left\|S^{1/2}G^\pm_\e f\right\|^2 &\leq  c_1\left|\Re \left\langle
G^\pm_\e f,(H^\pm_\varepsilon -\lambda  \pm \rmi \mu)  G^\pm_\e
f\right\rangle\right|+  
\frac{1}{\e}\left|\Im \left\langle G^\pm_\e f,(H^\pm_\varepsilon -\lambda  \pm
\rmi \mu) G^\pm_\e 
f\right\rangle\right| 
\\
\label{e:estG}
&\leq \max\left(c_1,\frac{1}{\e}\right)\left|F_\e^\pm\right|.
\end{align}
Hence up to a smaller $\e_0>0$, we obtain $\left\|S^{1/2}G^\pm_\e
f\right\|^2\leq  \left|F_\e^\pm\right|/ \e$ for all $\varepsilon
\in(0, \varepsilon_0]$. Moreover, if $f\in\Dc(S^{-1/2})$, we obtain
\begin{equation*}
	\left|F_\e^\pm\right|\leq \left\|S^{-1/2} f\right\|
	\left\|S^{1/2}G^\pm_\e f\right\|\leq  
	\left\|S^{-1/2} f\right\|
	\frac{ \sqrt{\left|F_\e^\pm\right|}}{\sqrt{\e}}
\end{equation*}
and deduce
\begin{equation}\label{e:estFe}
	\left|F_\e^\pm\right|\leq  
	\frac{1}{\e}\left\|S^{-1/2} f\right\|^2, \mbox{ for all } \varepsilon
\in(0, \varepsilon_0] .
\end{equation}
We now show that $G_\varepsilon^\pm\in \Cc^1(A)$. First note that
$G_\varepsilon^\pm$ is a bijection from $\Hr$ onto $\Dr$. Then by taking
the adjoint, it is also a bijection from $\Dr^*$ onto $\Hc$. Remember 
now that $W_t$ stabilizes $\Gr$ and $\Gr^*$. By the resolvent equality
in $\Bc(\Hr)$, we have: 
\begin{eqnarray*}
[G_\varepsilon^\pm , W_t]= -
\underbrace{ G_\varepsilon^\pm}_{\Hr\longleftarrow \Gr^*}\quad 
\underbrace{\big[H^\pm \pm\rmi \e[ H^\pm,\rmi A], W_t\big]}_{\Gr^*
  \longleftarrow \Gr}\quad 
\underbrace{G_\varepsilon^\pm}_{\Gr \longleftarrow \Hr} 
\end{eqnarray*}  
Let now take the derivative in $0$. Since $H^{\pm}$ and $[ H^\pm,\rmi A]$ are in
$\Cc^1(A; \Gr, \Gr^*)$ (the former being in $\Cc^2(A; \Gr, \Gr^*)$),  
the right hand side has a strong limit for all element in
$\Hr$. Hence, $G_\varepsilon^\pm \in \Cc^1(A; \Hr, \Hr)$ which is the same as
$G_\varepsilon^\pm \in \Cc^1(A)$, see
\cite{AmreinBoutetdeMonvelGeorgescu}[Theorem 6.3.4 a.]. 
Easily, it follows that $G_\varepsilon^\pm \Dc(A)\subset \Dc(A)\cap \Dr$ and
one can safely expand  the commutator in the next computation. Take
$f\in\Dc(A)$.  
\begin{align*}
\frac{d}{d\e}F_\e^\pm&=
\left\langle f, \frac{d}{d\e}G^\pm_\e f\right\rangle
=\pm\rmi \left\langle G^\mp_\e f, [ H^\pm,\rmi A]_\circ G^\pm_\e f\right\rangle\\
&=\pm \left\langle  G^\mp_\e f,Af\right\rangle \mp \left\langle
Af, G^\mp_\e f\right\rangle - \varepsilon  
\left\langle G^\mp_\e f, \big[[H,iA]_\circ, iA\big] G^\pm_\varepsilon f
\right\rangle. 
\end{align*}
Here the last commutator in taken in the form sense. 
Now use three times \eqref{e:estG} and the bound
\eqref{Eq:SecondCommBound}, which is uniform in $p\in \Ec$, then integrate 
to obtain
\begin{equation}\label{e:inegdiff}
\left|F_\e^\pm - F_{\e'}^\pm\right|\leq \int_\e^{\e'}\left\{2 \frac{
  \sqrt{|F_s^\pm|}} 
{\sqrt{s}}\left\|S^{-1/2}A f\right\| 
+ C \left|F_s^\pm\right|\right\}\;ds, \mbox{ for }
0<\varepsilon\leq\varepsilon'\leq \varepsilon_0
\end{equation}
and for all $f\in\Dc(S^{-1/2}A)\cap\Dc(A)$. 

We give a first
estimation. Using \eqref{e:estFe} and  the Gronwall lemma, see  
\cite{AmreinBoutetdeMonvelGeorgescu}[Lemma 7.A.1] 
with $\theta=1/2$ or \cite{Oguntuase}[Lemma 2.6] with $p=1/2$,
we infer there are some constants $C, C', C'', C'''$, independent of 
$\e\in(0,\e_0]$, $\lambda\geq 0$, $\mu>0$ and of $p\in \Ec$, so that
\begin{align}
\nonumber
\left|F_\e^\pm\right|&\leq e^{C(\e-\e_0)}
\left(\left|F_{\e_0}^\pm\right|^{1/2} 
+ \int_\e^{\e_0}\left\{\frac{  1}{\sqrt{\eta}}
e^{-\frac{1}{2}C(\eta-\e_0)}\right\}d\eta\,\,
\left\|S^{-1/2}A f\right\|\right)^2
\\ \nonumber
&\leq C''\left(\left|F_{\e_0}^\pm\right| 
+\left(\sqrt{\e}-\sqrt{\e_0}\right)^2
\left\|S^{-1/2}A f\right\|^2\right)
\\ \label{e:last}
&\leq C''\left(\frac{1}{\e_0}
\left\|S^{-1/2} f\right\|^2 
+\left(\sqrt{\e}-\sqrt{\e_0}\right)^2
\left\|S^{-1/2}A f\right\|^2\right)
\leq C''' \|f\|^2_{\tilde \Sr^*}
\end{align}
for $f\in\Dc(S^{-1/2})\cap\Dc(S^{-1/2}A)\cap\Dc(A)$
and where $\tilde \Sr^*$ is the completion of 
$\Dc(A|_{\Sr^*})$ under the norm $\|f\|^2_{\tilde \Sr^*}:=
\left\|S^{-1/2} f\right\|^2  +  \left\|S^{-1/2}A f\right\|^2$. 
Here one notices that the norm is well defined for elements of
$\Dc(A|_{\Sr^*})$ by taking in account \eqref{e:dom}. 
We now
plug this back in 
\eqref{e:inegdiff}. Since the inverse of the square root is integrable
around $0$, we find $C''''$ with the same independence so that 
\begin{align*}
\left|F_\e^\pm-F_{\e'}^\pm\right| 
\leq \int_\e^{\e'}\left\{2 \frac{\sqrt{C'''}}{\sqrt{s}} 
+ CC''' \right\}ds\,\, \|f\|_{\tilde \Sr^*}^2= C''''
\big(\sqrt{\varepsilon'}-\sqrt{\varepsilon}\,\big)  \|f\|_{\tilde \Sr^*}^2.
\end{align*}
Then, $\{F_\e^\pm\}_{\e\in(0, \varepsilon_0]}$ is a Cauchy sequence. We denote by 
$F_{0^+}^\pm$ the limit, as $\varepsilon$ goes to $0$. It remains to
notice that $F_{0^+}^\pm= F_0^{\pm}$. Indeed, using \eqref{e:G1} and
\eqref{e:mourrestrict0}, one has the stronger fact that
\begin{eqnarray*} 
\|G_0^\pm - G_\e^\pm\|_{\Bc(\Hr)}\leq \varepsilon \|G^\pm_\e\|_{\Bc(\Hr)}\cdot \| 
 [H^{\pm}(p),iA] (H^{\pm}(p)-\lambda\pm \rmi\mu)^{-1} 
\|_{\Bc(\Hr)} \leq \frac{c\varepsilon}{\mu^2}. 
\end{eqnarray*}
This gives us \eqref{e:LAP}. \qed

For the convenience of the reader, we give a proof of the following
well known fact:
\begin{lemma}[Numerical Range Theorem]\label{l:NRT} Let $H$ be a
  closed operator. Suppose that $\Dr:=\Dc(H)=\Dc(H^*)$. The numerical
  range of $H$ is defined by $\Nc:=\{ \langle f, Hf \rangle$ with $
  f\in \Dr $ and $ \|f\|=1\}$. We have that $\sigma(H) \subset
  \overline{\Nc}$, the closure of $\Nc$. Moreover, if $\lambda \notin
  \sigma(H)$, then $\|(H-\lambda)^{-1}\|\leq 1/d(\lambda, \Nc)$.
\end{lemma}
\proof Let $\lambda \notin \overline{\Nc}$. There is $c:=d(\lambda,
\Nc)>0$, such that
$|\langle f, Hf \rangle - \lambda|\geq c$. Then,
\begin{eqnarray*}
\|(H-\lambda)f \|\geq c \|f\|, \quad \|(H^*-\overline\lambda)f \|\geq c \|f\|,
\end{eqnarray*}
for all $f\in \Dr$ and $\|f\|=1$. From the second part, we get  the range of
$(H-\lambda)$ is dense. Then, since $H$ is closed, the first part
gives that the range of $(H-\lambda)$ is closed. Hence, using again the
first inequality, $H-\lambda$ is bijective. The open mapping theorem
concludes. \qed

\section{Application to non-relativistic dispersive
  Hamiltonians}\label{s:nonrela} 

In this section, we give  an immediate application to the theory
exposed in Appendix \ref{s:WeakMourreTheory}. We do not discuss the
uniformity with respect to the external parameter. The latter would be
used in the heart of our approach, see Section \ref{s:positive}. We
discuss shortly the Helmholtz equation, see \cite{BenamouCastella,
  BenamouLafitte, Wang2, WangZhang}. In \cite{Royer}, one studies the
size of the resolvent of  
\begin{eqnarray*}
H_h := -h^2 \Delta +V_1(Q) -ih V_2(Q), \mbox{ as } h\rightarrow 0.
\end{eqnarray*}
This operator models accurately the propagation of the
electromagnetic field of a laser in material medium. The important
improvement between \cite{Royer} and the previous ones, is 
that he allows $V_2$ to be a smooth function tending to $0$ without
any assumption on the size of $\|V_2\|_\infty$. Note he supposes the
coefficients are smooth as some pseudo-differential calculus is used
to applied the non self-adjoint Mourre theory he develops. Then, he
discusses trapping conditions in the spirit of  \cite{Wang2}. Here, we
will stick to the quantum case and choose $h=-1$. To simplify the
presentation and expose some key ideas of Section \ref{s:positive}, we
focus on $L^2(\R^n; \C)$, with $n\geq 3$. For dimensions $1$ and $2$,
one needs to adapt the first part of (H2) and the weights  in
\eqref{e:nonrelamain}.   

\begin{theorem}\label{t:nonrelamain}
Suppose that $V_1, V_2\in L^1_{\rm loc}(\R^n;\R)$ satisfy:
\begin{enumerate}
\item[(H0)] $V_i$ are $\Delta$-operator bounded with a relative bound
  $a<1$, for  $i\in\{1,2\}$. 
\item[(H1)] $\nabla V_i$,  $ Q\cdot \nabla V_i(Q)$ are in
  $\Bc(\Hr^2(\R^n);L^2(\R^n))$ and $\langle Q\rangle (Q\cdot \nabla
  V_i)^2(Q)$ is bounded, for  $i\in\{1,2\}$.   
\item[(H2)] There are $c_1\in [0,2)$ and $\displaystyle c_1'\in
  \big[0, 4(2-c_1)/(n-2)^2\big)$ such that 
\begin{eqnarray*}
W_{V_1}(x):= x\cdot
  (\nabla V_1)(x) + c_{1} V_1(x) \leq \frac{c_{1}'}{|x|^2}, \mbox{ for
  all } x\in \R^n. 
\end{eqnarray*} 
and 
\begin{eqnarray*}
V_2(x)\geq 0 \mbox{ and } -c_1 x\cdot
  (\nabla V_2)(x)\geq  0, \mbox{ for all } x\in \R^n.
\end{eqnarray*} 
\end{enumerate}
On $\Cc^\infty_c(\R^n)$, we define $H:= -\Delta + V(Q), \mbox{ where }
V:=V_1 + i V_2$. The closure of $H$ defines a dispersive closed operator
with domain $\Hr^2(\R^n)$. We keep denoting it with $H$. Its spectrum
included in the upper half-plane. Moreover, $H$ has no eigenvalue in $[0,\infty)$ and
\begin{eqnarray}\label{e:nonrelamain}
\sup_{\lambda\in [0,\infty), \, \mu>0} \big\|\, |Q|^{-1}(H -\lambda +
  \rmi \mu )^{-1} |Q|^{-1}\big\| <\infty.	  
\end{eqnarray}	
If $c_1=0$, $H$ has no eigenvalue in $\R$ and \eqref{e:nonrelamain} holds true for $\lambda \in \R$. 
\end{theorem}
The quantity $W_{V_1}$ is called the \emph{virial} of $V_1$. 
For $h$ fixed and for a compact $\Ic$ included in
$(0,\infty)$, \cite{Royer} shows some estimates of the resolvent above
$\Ic$. Here we deal with the threshold $0$ and with high energy
estimates. On the other hand, as he avoids the threshold, he reaches
some very sharp weights.  As mentioned above, one can improve the
weights $|Q|$ to some extend by the use of Besov spaces, see
\cite{Richard}. In \cite{Royer} one makes an hypothesis on the sign of $V_2$ but not 
on the one of $x\cdot (\nabla V_2)(x)$. Note that if one supposes $c_1=0$, we are also in this situation. We take the opportunity to point out \cite{Wang},
where one discusses the presence of possible eigenvalues in $0$ for
non self-adjoint problems.  

\begin{remark}\label{r:virial} 
Taking $V_2=0$, we can compare the results with
\cite{FournaisSkibsted,Richard}.  
In \cite{FournaisSkibsted}, one uses in a crucial way that
$W_{V_1}(x)\leq - c \langle x\rangle^\alpha$ in a neighborhood of
infinity, for some $\alpha, c>0$. In \cite{Richard}, one remarks that
the condition $W_{V_1}(x)\leq 0$ is enough to obtain the
estimate. Here we mention that the condition (H2) is sufficient.  Note
this example is not explicitly discussed in \cite{Richard} but is
covered by his abstract approach. In \cite{BoutetKazantsevaMantoiu},
for the special case $c_1=0$,   one uses extensively the condition
(H2). This implies \eqref{e:nonrelamain} for $\lambda\in \R$. 
\end{remark} 

\begin{remark}\label{r:Herbst}
Unlike in \cite{Royer}, we stress that $V$ is \emph{not} supposed to
be a relatively compact perturbation of $H$ and that the essential
spectrum of $H$ can be different of $[0, \infty)$. 
In \cite{Herbst}, see also \cite{BoutetKazantsevaMantoiu}, one 
studies $V_2=0$ and $V_1(x):= v(x/|x|)$, with $v\in
\Cc^\infty(S^{n-1})$. We improve the weights of 
\cite{Herbst}[Theorem 3.2] from $\langle Q\rangle $ to $|Q|$. We
can also give a non-self-adjoint version. Consider 
$V_1$ satisfying (H1) and being relatively compact with respect to
$\Delta$ and $V_2(x):= v(x/|x|)$, where $v\in \Cc^0(S^{n-1})$, non-negative.
If $v^{-1}(0)$ is non-empty, one shows $[0,\infty)$ is included in the
  essential spectrum of $H$ by using some Weyl sequences.  
\end{remark} 

\proof[Proof of Theorem \ref{t:nonrelamain}.] Using (H0) and 
adapting the proof of Kato-Rellich, e.g.,
\cite{ReedSimon}[Theorem X.12],
one obtains easily $\Dc(H)=\Dc(H^*)=\Hr^2(\R^n)$. Let
$S:=c_s(-\Delta)^{1/2}$, with  $c_s:=2-c_{1}- (n-2)^2
c_{1}'/4>0$. Set $\Sr:=\dot\Hr^1(\R^n)$, the homogeneous Sobolev
space of order $1$, i.e., the completion of  $\Hr^1(\R^n)$ under the
norm $\|f\|_\Sr:= \| S^{1/2} f\|^2$. Consider the strongly continuous
one-parameter unitary group $\{W_t\}_{t\in \R}$ acting by: $(W_t
f)(x)= e^{nt/2} f(e^tx)$, for all  $f\in L^2(\R^3)$. This is the
$C_0$-group of dilatation. By interpolation and 
duality, one derives  $W_t \Sr \subset \Sr \mbox{ and } W_t \Hr^s(\R^3) \subset
\Hr^s(\R^n)$, for all $s\in \R$. Consider now its generator $A$ in $L^2(\R^n)$. 
By the Nelson lemma, it is essentially self-adjoint on
$\Cc^\infty_c(\R^n)$ and acts as follows: $A=(P\cdot Q+ Q\cdot P)/2$ on
$\Cc^\infty_c(\R^n)$. By computing on $\Cc^\infty_c(\R^n)$ in the form
sense, we obtain that  
\begin{eqnarray}\label{e:nonrelamain1}
[\Re(H), iA] - c_{1} \Re(H) = -(2-c_{1}) \Delta - W_{V_1}\geq S,
\end{eqnarray}
here we used the Hardy inequality for the last step. Furthermore,
$\Im(H)=V_2(Q)\geq 0$, 
\begin{eqnarray}\label{e:nonrelamain2}
[\Im(H), iA]= -Q\cdot\nabla(V_2)(Q),
\end{eqnarray}
and also
\begin{eqnarray}\label{e:nonrelamain3}
[[H, iA], iA] = -4 \Delta + (Q\cdot \nabla V)^2(Q).
\end{eqnarray}
Since $W_t$ stabilizes $\Gr:=\Hr^1$ and as \eqref{e:nonrelamain1},
\eqref{e:nonrelamain2} and \eqref{e:nonrelamain3} extend to bounded
operators from $\Hr^1$ into $\Hr^{-1}$, we infer that $H$ and $H^*$
are in $\Cc^2(A; \Hr^1, \Hr^{-1})$ and also  \eqref{e:mourrestrict1}
and \eqref{e:mourrestrict2}. 
Now since $\Cc^\infty_c(\R^3)$ is a core for $H$, $H^*$ and $A$,
\eqref{e:nonrelamain1} and \eqref{e:nonrelamain2} give
\eqref{e:mourrestrict0}, with notation $H^+=H$ and $H^-=H^*$. In
addition \eqref{Eq:SecondCommBound} follows from the Hardy inequality
and (H1), as $\|(Q\cdot \nabla)^2 v(Q) f\|^2\leq c \|\, 
|Q|(Q\cdot \nabla)^2 v(Q)\|^2 \|S f\|^2$. Therefore, we can apply
Theorem \ref{t:mourrestrict} and derive the weight $|Q|$ by the Hardy
inequality. \qed

Finally, we recall the Hardy inequality. Take $E$ a finite
dimensional vector space. One has:
\begin{eqnarray}\label{e:Hardy}
\left(\frac{n-2}{2}\right)^2\int_{\R^n} \left| \frac{1}{|x|}
f(x)\right|^2\, dx \leq \, \big|\langle f, -\Delta f \rangle\big|,
\mbox{ where } n\geq 3 \mbox{ and } f\in \Cc^\infty_c(\R^n; E). 
\end{eqnarray}
 
\bibliographystyle{amsplain}

\providecommand{\bysame}{\leavevmode\hbox to3em{\hrulefill}\thinspace}
\providecommand{\MR}{\relax\ifhmode\unskip\space\fi MR }
\providecommand{\MRhref}[2]{%
 \href{http://www.ams.org/mathscinet-getitem?mr=#1}{#2}
}
\providecommand{\href}[2]{#2}

\end{document}